\documentclass[11pt,reqno,a4]{amsart}

\usepackage{a4}
\usepackage{graphicx}
\usepackage{amsmath,amssymb,color,bbm,ifthen,enumerate}
\usepackage{hyperref}
\usepackage[latin1]{inputenc}
\usepackage[numbers]{natbib}

\newtheorem{theorem}{Theorem}
\newtheorem{lemma}[theorem]{Lemma}

\newtheorem{proposition}[theorem]{Proposition}
\newtheorem{corollary}[theorem]{Corollary}
\newtheorem{remark}[theorem]{Remark}
\makeatletter \@addtoreset{equation}{section} \makeatother
\renewcommand{\theequation}{\arabic{section}.\arabic{equation}}
  \advance\oddsidemargin -1.4cm
  \advance\evensidemargin -1.4cm

\def\R{{\mathbb R}}
\def\Q{{\mathbb Q}}

\def\dim{{\rm dim}_{_{\rm H}}}

\def\dimr{{\rm dim}^\rho_{_{\rm H}}}
\def\E{{\mathbb E}}
\def\N{{\mathbb N}}
\def\P{{\mathbb P}}
\def\Z{{\mathbb Z}}
\def\l{{\langle}}
\def\r{\rangle}
\def\rme{\mathrm{e}}

\def\1{\mathbbm{1}}
\def\a{\alpha}

\def\de{\delta}
\def\ga{\gamma}

\def\ep{\varepsilon}

\def\K{\mathcal{K}}
\def\cqfd{$\square$}

\title[Linear Fractional Stable Sheet]{Linear Fractional Stable Sheets:
wavelet expansion and sample path properties}

\author{Antoine Ayache}
\address{UMR CNRS 8524, Laboratoire Paul Painlev\'e, B\^at. M2,
  Universit\'e Lille 1, 59655 Villeneuve d'Ascq Cedex, FRANCE}
\email{Antoine.Ayache@math.univ-lille1.fr}

\author{François Roueff}
\address{TELECOM ParisTech, CNRS LTCI, Institut Télécom, 46, rue Barrault, 75634 Paris Cédex 13, France.}
\email{roueff@tsi.enst.fr}

\author{Yimin Xiao}
\thanks{Yimin Xiao's research is partially supported by the NSF grant DMS-0706728.}

\address{Department of Statistics and Probability, Michigan State University,
East Lansing, MI 48824, U.S.A.}
\email{xiao@stt.msu.edu}



\date{June 9, 2008}

\begin{document}

\maketitle

\begin{center}
{\it UMR CNRS 8524, TELECOM ParisTech and Michigan State University}\\
\end{center}
\renewcommand{\thefootnote}{}
\footnote{\textit{Corresponding author}: Antoine Ayache
(\email{Antoine.Ayache@math.univ-lille1.fr})  }
\renewcommand{\thefootnote}{\arabic{footnote}}

\begin{abstract}

In this paper we give a detailed description of the random wavelet
series representation of real-valued linear fractional stable
sheet introduced in~\cite{ayache:roueff:xiao:2007a}. By using this
representation, in the case where the sample paths are continuous, an
anisotropic uniform and quasi-optimal modulus of continuity of these
paths is obtained  as well as an upper bound for their behavior at
infinity and around the coordinate axes. The Hausdorff dimensions of the
range and graph of these stable random fields are then derived.

\end{abstract}

\section{Introduction and main results}

Let $0 < \alpha < 2$ and $H = (H_1,
\ldots, H_N) \in (0, 1)^N$ be given. We define an $\a$-stable random
field $X_0 = \{X_0(t), t \in \R^N\}$ with values in $\R$ by
\begin{equation}\label{Eq:Repstable}
X_0(t) = \int_{\R^N} h_{_H} (t, s)\, Z_\a(ds),
\end{equation}
where $Z_\a$ is a strictly $\alpha$-stable random
measure on $\R^N$ with Lebesgue measure as its control measure and
$\beta(s)$ as its skewness intensity. That is, for every Lebesgue
measurable set $A \subseteq \R^N$ with Lebesgue measure $\lambda_N(A)
< \infty$, $Z_\a(A)$ is a strictly $\alpha$-stable random variable
with scale parameter $\lambda_N(A)^{1/\a}$ and skewness parameter
$(1/\lambda_N(A))\int_A \beta(s) ds.$  If $\beta(s) \equiv 0$, then
$Z_\a$ is a symmetric $\alpha$-stable random measure on $\R^N$. We
refer to \cite[Chapter 3]{ST94} for more
information on stable random measures and their integrals. Also in
(\ref{Eq:Repstable}),
\begin{equation}\label{Eq:functh}
h_{_H} (t, s) = \kappa\, \prod_{\ell=1}^N
\Big\{(t_\ell-s_\ell)_+^{H_\ell- \frac 1 \a} - (-s_\ell)_+^{H_\ell-
\frac 1 \a} \Big\}\;,
\end{equation}
where $\kappa> 0$ is a normalizing constant such that the scale
parameter of $X_0(1)$, denoted by $\|X_0(1)\|_\a$, equals 1,
$t_{+} = \max\{t, 0\}$ and $0^0=1$.
Observe that, if $H_1 = \cdots = H_N = \frac 1 {\a}$, $X_0$ is the ordinary
stable sheet studied in \cite{ehm}. In general,
the random field $X_0$ is called a linear fractional $\a$-stable
sheet defined on $\R^N$ (or $(N,1)$-LFSS for brevity) in $\R$ with
index $H$. LFSS is an extension  of both
linear fractional stable motion (LFSM), which corresponds to the
case where $N=1$, and ordinary fractional Brownian sheet (FBS)
which corresponds to $\a=2$, that is, to replacing the stable
measure in~(\ref{Eq:Repstable}) by a Gaussian random measure.

We will also consider $(N,d)$-LFSS, with $d>1$, that is
a linear fractional $\a$-stable sheet defined on $\R^N$ and taking
its values in $\R^d$. The $(N, d)$-LFSS that we consider is the
stable field $X = \{X(t), t \in \R^N\}$  defined by
\begin{equation}\label{Def:X}
X(t) = \big(X_1(t), \ldots, X_d(t)\big),\qquad \forall t \in
\R^N \; ,
\end{equation}
where $X_1, \ldots, X_d$ are $d$ independent copies of $X_0$.
It is easy to verify by using the representation
(\ref{Eq:Repstable}) that  $X$ satisfies the following scaling
property: For any $N\times N$ diagonal matrix $A = (a_{ij})$ with
$a_{ii} = a_i > 0$ for all $1 \le i \le N$ and $a_{ij}=0$ if $i\ne
j$, we have
\begin{equation}\label{Eq:OSS}
\big\{ X(A t),\, t \in \R^N \big\} \stackrel{d}{=} \bigg\{
\bigg(\prod_{j=1}^N a_j^{H_j} \bigg)\,X(t),\ t \in \R^N \bigg\},
\end{equation}
where $\stackrel{d}{=}$ denotes the equality in the sense of finite
dimensional distributions, provided that the skewness
intensity satisfies $\beta(As)=\beta(s)$ for almost every $s\in\R^N$.
Relation~(\ref{Eq:OSS}) means that the $(N, d)$-LFSS $X$ is an
operator-self-similar [or operator-scaling] random field in the
\emph{time variable} (see ~\cite{BMS07, Xiao07}). When the indices
$H_1, \ldots, H_N$ are not the same, $X$ has different scaling behavior
along different directions.  This anisotropic nature of $X$
makes it a potential model for various spatial objects, as is
already the case for anisotropic Gaussian fields (\cite{bensoetal:2006}
and \cite{bonami:estrade:2003}). We also mention that one can construct
$(N, d)$-stable random fields which are self-similar in the
\emph{space variables} in the sense of  ~\cite{LX94, mason:xiao:2001}.
This will not be discussed in this paper.

Similarly to LFSM and FBS, see for
instance \cite{ayache2004,AX05,KoMae91,Maejima83,takashima:1989,Xiao06},
there are close connections between sample path
properties of LFSS and its parameters $H$ and $\alpha$. In this
article we study some of these connections. In all the
remainder of this paper we assume that the sample paths of $X$ are continuous,
\textit{i.e.} $\min(H_1,\dots,H_N)>1/\alpha$.  For convenience
we even assume that
\begin{equation}
\label{Eq:H}
1/\a < H_1 \leq \dots \leq H_N \; .
\end{equation}
Of course, there is no loss of generality in the arbitrary
ordering of $H_1,\dots,H_N$. Observe that, 
since $H\in(0,1)^N$, condition (\ref{Eq:H})
implies $\a>1$.

Let us now state our main results.  Theorem \ref{thm:mod-cont} below is an 
improved version of Theorems~1.2 and~1.3
in~\cite{ayache:roueff:xiao:2007a}. Relation~(\ref{eq:mod-cont}) provides
a sharp upper bound for the uniform modulus of continuity of LFSS, while
Relation~(\ref{eq:ub-infty}) gives an upper bound for its asymptotic
behavior at infinity and around the coordinate axes.

\begin{theorem}
\label{thm:mod-cont}
Let $\Omega_{0}^*$ be the event of probability $1$ that will be introduced
in Corollary~\ref{cor:ub}. Then for every compact
set ${\mathcal K}\subseteq \R^N$, all
$\omega\in\Omega_{0}^*$ and any arbitrarily small $\eta>0$, one has
\begin{equation}
\label{eq:mod-cont}
\sup_{s,t\in {\mathcal K}}\frac{|X_0(s,\omega)-X_0(t,\omega)|}{\sum_{j=1}^N
|s_j-t_j|^{H_j-1/\alpha}\big(1+\big|\log |s_j-t_j|\big| \big)^{2/\a+\eta}}<\infty
\end{equation}
and
\begin{equation}
\label{eq:ub-infty}
\sup_{t\in\R^N}\frac{|X_0(t,\omega)|}{\prod_{j=1}^N
|t_j|^{H_j}(1+\big|\log|t_j| \big|)^{1/\alpha+\eta}}<\infty.
\end{equation}
\end{theorem}

The following result can be viewed as an inverse
of~(\ref{eq:mod-cont}) in Theorem~\ref{thm:mod-cont}.
\begin{theorem}
\label{thm:irreg} Let $\Omega_{3}^*$ be the event of probability $1$
that will be introduced in Lemma \ref{lem:liminf-G}. Then for all
$\omega\in\Omega_3^*$, all vectors $\widehat{u}_n\in\R^{N-1}$ with
non-vanishing coordinates, any $n=1,\ldots, N$ and any real numbers
$y_1<y_2$ and $\epsilon>0$, one has
\begin{equation}
\label{eq:irreg} \sup_{s_n,t_n\in
[y_1,y_2]}\frac{|X_0(s_n,\widehat{u}_n,\omega)-X_0(t_n,\widehat{u}_n,\omega)|}
{|s_n-t_n|^{H_n-1/\alpha}\big(1+\big|\log |s_{n}-t_{n}| \big| \big)^{-1/\a-\epsilon}}
=\infty,
\end{equation}
where, for every real $x_n$, we have set
$(x_n,\widehat{u}_n)=(u_1,\ldots,u_{n-1},x_n,u_{n+1},\ldots,
u_N)$.
\end{theorem}
Observe that Theorems \ref{thm:mod-cont} and \ref{thm:irreg} have
already been obtained by Takashima \cite{takashima:1989}
in the particular case of LFSM (\textit{i.e.}, $N=1$).
However, the proofs given by this author can hardly be adapted to LFSS.
To establish the above theorems we introduce
a wavelet series representation of
$X_0$ and use wavelet methods which are, more or less, inspired from~\cite{AX05}.
It is also worth noticing that  the event $\Omega_{3}^*$ in Theorem
\ref{thm:irreg} does not depend on $\widehat{u}_n$. This is why the latter
theorem cannot be obtained by simply using the fact that LFSS is an LFSM
of Hurst parameter $H_n$ along the direction of the $n$-th axis.

The next theorem gives the Hausdorff dimensions of the range
$$
X\big([0, 1]^N\big) = \big\{X(t):\, t \in [0, 1]^N \big \}
$$
and the graph
$$
{\rm Gr}X\big([0, 1]^N\big) = \big\{(t, X(t)):\, t \in
[0, 1]^N\big\}
$$
of an $(N, d)$-LFSS $X$. We refer to~\cite{falconer90} for the
definition and basic properties of Hausdorff dimension.

The following result extends Theorem 4 in \cite{AX05}
to the linear fractional stable sheets. Unlike the fractional
Brownian sheet case, we remark that the Hausdorff dimensions of
$X\big([0, 1]^N\big)$ and ${\rm Gr}X\big([0, 1]^N\big)$ are not
determined by the uniform H\"older exponent of $X$ on $[0, 1]^N$.

\begin{theorem} \label{Th:dim}
Let the assumption (\ref{Eq:H}) hold. Then, with probability 1,
\begin{equation}\label{Eq:range}
\dim X\big([0, 1]^N\big) = \min \bigg\{ d; \quad \sum_{\ell=1}^N
\frac1 {H_\ell}  \bigg\}
\end{equation}
and
\begin{equation}\label{Eq:graph}
\begin{split}
&\ \dim {\rm Gr}X\big([0, 1]^N\big) = \min \bigg\{ \sum_{\ell=1}^k
\frac{H_k} {H_\ell} + N-k +
(1 - H_k) d,\ 1 \le k \le N; \ \sum_{\ell=1}^N \frac1 {H_\ell}\bigg\}  \\
&\  \quad = \left\{ \begin{array}{ll}
\sum_{\ell=1}^N \frac1 {H_\ell}  &\hbox{ if }\
\sum_{\ell=1}^N \frac 1 {H_\ell} \le d, \\
\sum_{\ell=1}^k \frac{H_k} {H_\ell} + N-k + (1 - H_k) d  &\text {
if } \ \sum_{\ell=1}^{k-1} \frac 1 {H_\ell} \le d <
\sum_{\ell=1}^{k} \frac 1 {H_\ell}\, ,
\end{array}\right.
\end{split}
\end{equation}
where $\sum_{\ell=1}^{0} \frac 1 {H_\ell} := 0$.
\end{theorem}

\begin{remark}
{\rm The second equality in (\ref{Eq:graph}) can be verified by
using (\ref{Eq:H}) and some elementary
computation; see \cite{AX05}.}
\end{remark}

In light of Theorem \ref{Th:dim} it is a natural question to consider
the Hausdorff dimensions of the image $X(E)$ and graph ${\rm Gr}X(E)$,
where $E$ is an arbitrary Borel set in $\R^N$. As shown by Wu and Xiao
~\cite{wu:xiao:2007} for fractional Brownian sheets, due to the anisotropic
nature of $X$, the Hausdorff dimension of $E$ and the index $H$ alone
are not enough to determine $\dim X(E)$. By combining the methods in
Wu and Xiao \cite{wu:xiao:2007} and Xiao \cite{Xiao07} with the moment
argument in this paper we determine $\dim X(E)$ for every nonrandom Borel
set $E \subseteq (0, \infty)^N$; see Theorem \ref{Th:dimXE}.

We end the Introduction with some notation. Throughout this paper,
the underlying parameter spaces are $\R^N$, $\R^N_+ = [0,
\infty)^N$ or $\Z^N$.  A typical parameter, $t\in\R^N$ is written as
$t = (t_1, \ldots, t_N)$ or $t = \l t_j\r$ whichever is more convenient.
For any $s, t \in \R^N$ such that $s_j < t_j$ ($j = 1,
\ldots, N$), the set $[s, t] = \prod^N_{j=1}\, [s_j, t_j]$ is called a
closed interval (or a rectangle). Open or half-open intervals can be
defined analogously. We will use capital letters $C, C_1, C_2, \ldots$ to
denote positive and finite random variables and use $c, c_1, c_2, \ldots$
to denote unspecified positive and finite constants. Moreover, $C$ and $c$
may not be the same in each occurrence.

\medskip
\noindent{\bf Acknowledgment}\ This paper was finished while the third author
(Y. Xiao) was visiting the Statistical \& Applied Mathematical Sciences Institute
(SAMSI). He thanks Professor Jim Berger and the staff of SAMSI for their
support and the good working conditions.

\section{Wavelet expansion of LFSS}

The goal of this section is to give a detailed description of the
wavelet representations of LFSS $X_0$. First we need to introduce some
notation that will be extensively used in all the sequel.
\begin{itemize}
\item[(i)] The real-valued function $\psi$ denotes a well chosen
compactly supported Daubechies wavelet (see~\cite{daubechies:1992,meyer:1992}).
Contrary to the Gaussian case
the fact that $\psi$ is compactly supported will play a crucial role
in the proof of Theorem  \ref{thm:irreg} (see the proof of Part
$(b)$ of Proposition \ref{prop:proper-G}).
\item[(ii)] For any $\ell=1,\dots,N$, the real-valued functions $\psi^{H_\ell}$ and
  $\psi^{-H_\ell}$ respectively denote the left-sided fractional primitive of order
  $H_\ell+1-1/\alpha$ and the right-sided fractional derivative of order
  $H_\ell+1-1/\alpha$ of $\psi$, which are respectively defined for all $x\in\R$ by
\begin{equation}
\label{eq:psi-frac}
\psi^{H_\ell} (x)=\int_{\R}(x-y)_{+}^{H_\ell-1/\alpha}\psi(y)\;dy
\mbox{ and }
\psi^{-H_\ell}(x)=\frac{d^2}{dx^2}\int_{\R}(y-x)_{+}^{1/\alpha-H_\ell}\psi (y)\, dy.
\end{equation}
Observe that the functions $\psi^{H_\ell}$ and $\psi^{-H_\ell}$  are
well-defined, continuously differentiable and well-localized
provided that $\psi$ has sufficiently many vanishing moments (and thus
is smooth enough).
By saying that a function $\phi: \R \to \R$ is
well-localized we mean that
\begin{equation}
\label{eq:wl-frac} \sup_{x\in\R} (1+|x|)^{2} \left\{ \left
|\phi(x)\right|+  \left |\phi'(x)\right| \right\} < \infty \;.
\end{equation}

\item[(iii)] $\{\epsilon_{j,k},\, (j,k)\in\Z^N\times\Z^N\}$ will denote the sequence
  of random variables defined as
\begin{equation}
\label{eq:def-epsilon} \epsilon_{j,k}=\int_{\R^N}\prod_{\ell=1}^N
\left\{2^{j_\ell/\alpha}\psi(2^{j_\ell}s_\ell-k_\ell)\right\}\,Z_\alpha(ds)\, .
\end{equation}
They are strictly $\a$-stable 
random variables all with the same scale parameter
$$
\|\epsilon_{j,k}\|_\a= \left\{\int_{\R} |\psi(t)|^\alpha\,dt\right\}^{N/\alpha} 
$$
and skewness parameter
$$
\beta_{j,k}=\|\epsilon_{j,k}\|_\a^{-\a}\int_{\R^N}\prod_{\ell=1}^N
\left\{2^{j_\ell}\psi^{<\a>}(2^{j_\ell}s_\ell-k_\ell)\right\}\beta(s)\,ds\;,  
$$
where $x^{<\a>} = |x|^\a\hbox{sgn}(x)$ which is the number having the same sign as
$x$ and absolute value $|x|^\a$.
Moreover, if $L>0$ is a constant such that the support of $\psi$ is included in $[-L,L]$,
then for any integers $p>2L$, any $r\in\{0,\dots,p-1\}^N$ and $j\in\Z^N$,
$\{\epsilon_{j,r+kp};\;k\in\Z^N\}$ is a sequence of independent random
variables.
\end{itemize}

A consequence of the above properties of the sequence $\{\epsilon_{j,k},\,
(j,k)\in\Z^N\times\Z^N\}$ is the following.

\begin{corollary}
\label{cor:ub}
There exists an event $\Omega_0^*$ of probability 1 such that, for any $\eta>0$,
for all $\omega\in\Omega_0^*$
and all $j,k\in\Z^N\times\Z^N$,
$$
|\epsilon_{j,k}(\omega)|\leq C(\omega) \prod_{l=1}^N \left\{(1+|j_l|)^{1/\a+\eta}
(1+|k_l|)^{1/\a}\log^{1/\a+\eta}(2+|k_l|)\right\} \;,
$$
where $C$ is a finite positive random variable.
\end{corollary}
\noindent {\bf Proof.}\,
We apply Lemma~\ref{lem:up-ps-indep}.\cqfd

It is worth noticing that, for every $\ell=1,\ldots, N$, the functions
$\psi^{H_\ell}$ and $\psi^{-H_\ell}$ can be
defined equivalently to~(\ref{eq:psi-frac}), up to a multiplicative
constant, but in the Fourier domain by (see
\textit{e.g.}~\cite{samkoetal:1993})
\begin{equation}
\label{eq:ft-fracwv1} \widehat{\psi^{H_\ell}}(\xi)=\rme^{i\,\mbox{{\Small
sgn}}(\xi)(H_\ell-1/\alpha+1)
\frac{\pi}{2}}\frac{\widehat{\psi}(\xi)}{|\xi|^{H_\ell-1/\alpha+1}}
\end{equation}
 and
 \begin{equation}
\label{eq:ft-fracwv2} \widehat{\psi^{-H_\ell}}(\xi)=\rme^{i\,\mbox{{\Small
sgn}}(\xi)(H_\ell-1/\alpha+1) \frac{\pi}{2}}
|\xi|^{H_\ell-1/\alpha+1}\widehat{\psi}(\xi).
\end{equation}
It follows from Parseval's Formula, (\ref{eq:ft-fracwv1}),
(\ref{eq:ft-fracwv2})  and the orthonormality (in  $L^2(\R)$) of
the sequence $\{2^{j/2}\psi(2^j\cdot-k),\,j,k\in\Z\}$ that
$\psi^{H_\ell}$ and $\psi^{-H_\ell}$  satisfy, for all $(J,K)\in\Z^2$
and $(J',K')\in\Z^2$, up to a multiplicative constant,
\begin{equation}
\label{eq:biorth-psiH} \int_{\R}\psi^{H_\ell} (2^Jx-K)
\overline{\psi^{-H_\ell}(2^{J'}x-K')}\,dx=2^{-J}\delta (J,K;J',K'),
\end{equation}
where $\delta (J,K;J',K')=1$ when $(J,K)=(J',K')$ and $0$ otherwise.
By putting
together (\ref{eq:ft-fracwv1}), (\ref{eq:ft-fracwv2}) and the fact
that $\widehat{\psi}(\xi)=O(\xi^2)$ as $|\xi|\to0$, another useful
property is obtained: for every $\ell=1,\ldots, N$, the first
moment of the functions  $\psi^{H_\ell}$ and $\psi^{-H_\ell}$ vanish,
namely one has
\begin{equation}
\label{eq:moment-psiH}
\int_{\R}\psi^{H_\ell}(u)\,du=\int_{\R}\psi^{-H_\ell}(u)\,du=0.
\end{equation}

We are now in position to state the main results of this section.

\begin{proposition}
\label{prop:wav-holder}
Let $\Omega_{1}^*$ be the event of probability $1$ that will be introduced in
Lemma \ref{lem:ub-eps}.
For every $n\in\N$, $M>0$ and $t\in\R^N$ we set
\begin{equation}
\label{eq:def-Un}
U_{n,M}(t)=\sum_{(j,k)\in D_{n,M}^N}2^{-
\langle j,H \rangle }\epsilon_{j,k}\prod_{l=1}^N
\left\{\psi^{H_l}(2^{j_l}t_l-k_l)-\psi^{H_l}(-k_l)\right\} \;,
\end{equation}
where the random variables $\{\epsilon_{j,k},\, (j,k)\in
\Z^N\times\Z^N\}$ are defined by~(\ref{eq:def-epsilon}) and
\begin{equation}
\label{eq:def-D-Dn}
D_{n,M}^N=\Big\{(j,k)\in\Z^N\times\Z^N\,:\, \mbox{for all $l=1,\ldots,N$ }|j_l|
\le n \mbox{ and } |k_l|\le M2^{n+1}\Big\}\;.
\end{equation}
Then for every $\omega\in\Omega_{1}^*$ the functional sequence
$(U_{n,M}(\cdot,\omega))_{n\in \N}$ is a Cauchy sequence in
the H\"older space $C^\gamma (\K)$ for every $\gamma\in [0,H_1-1/\alpha)$
and compact set $\K\subseteq [-M,M]^N$. We denote its
limit by
$$
\sum_{(j,k)\in \Z^N\times \Z^N}2^{-
\langle j,H \rangle }\epsilon_{j,k}\prod_{l=1}^N
\left\{\psi^{H_l}(2^{j_l}t_l-k_l)-\psi^{H_l}(-k_l)\right\}.
$$
\end{proposition}

\begin{proposition}
\label{thm:wav-dec}
With probability 1, the following holds for  all $t \in \R^N$
\begin{equation}
\label{eq:wav-dec} X_0(t)=\sum_{(j,k)\in\Z^N\times\Z^N}2^{-
\langle j,H \rangle }\epsilon_{j,k}\prod_{l=1}^N
\left\{\psi^{H_l}(2^{j_l}t_l-k_l)-\psi^{H_l}(-k_l)\right\}.
\end{equation}
\end{proposition}

\begin{remark}
\label{rem:modif-lfss}
{\rm By the definition of $X_0$ and by Proposition~\ref{prop:wav-holder},
both sides of~(\ref{eq:wav-dec}) are continuous in
$t$ with probability 1. Hence, to prove Proposition~\ref{thm:wav-dec},
it is sufficient to show that
$$
\left\{\sum_{(j,k)\in \Z^N\times \Z^N}2^{-
\langle j,H \rangle }\epsilon_{j,k}\prod_{l=1}^N
\left\{\psi^{H_l}(2^{j_l}t_l-k_l)-\psi^{H_l}(-k_l)\right\},\,t\in\R^N\right\}
$$
is a modification of $X_0$.
This is a natural extension of the wavelet
series representations both of LFSM and FBS (see \cite{benassi:roux:2003,ayache2004,AX05})
and will be
called the \emph{random wavelet series representation} of LFSS.}
\end{remark}

Assume for a while that Proposition \ref{prop:wav-holder}
holds and let us prove Proposition
\ref{thm:wav-dec}.

\noindent {\bf Proof of Proposition \ref{thm:wav-dec}.}\, Let us
fix $l\in\{1,\ldots, N\}$. For any $(j_l,k_l)\in\Z\times\Z$ and $s_l\in\R$ we set
\begin{equation}
\label{eq:def-psijk}
\psi_{j_l,k_l}(s_l)= 2^{j_l/\alpha} \psi(2^{j_l}s_l-k_l).
\end{equation}
Since $\{\psi_{j,k},j\in\Z,k\in\Z\}$ is an unconditional basis
of $L^\alpha(\R)$ (see~\cite{meyer:1990}) and,
for every fixed $t_l\in\R$,  the function
$s_l \mapsto (t_l-s_l)_+^{H_l-1/\alpha}-(-s_l)_+^{H_n-1/\alpha}
\in L^\alpha(\R) \cap L^2(\R)$, one has
\begin{equation}
\label{eq:KtWavDecomp}
(t_l-s_l)_+^{H_l-1/\alpha}-(-s_l)_+^{H_l-1/\alpha}
=\sum_{j_l\in\Z}\sum_{k_l\in\Z} \kappa_{l,j,k}(t_l)
\psi_{j_l,k_l}(s_l) \;,
\end{equation}
where the convergence of the series in the RHS of~(\ref{eq:KtWavDecomp}),
as a function of $s_l$, holds in $L^\alpha(\R)$. Next by
using the H\"older inequality and the $L^2(\R)$ orthonormality
of the sequence $\big\{2^{j_l(1/2-1/\alpha)}\psi_{j_l,k_l},\,$
$j_l\in\Z,k_l\in\Z\big\}$, one can prove that
\begin{eqnarray}
\label{eq:def-kappa}
\nonumber
\kappa_{l,j_l,k_l}(t_l)&=& 2^{j_l(1-1/\alpha)}
\int_{\R} \{(t_l-s_l)_+^{H_l-1/\alpha}-(-s_l)_+^{H_l-1/\alpha}\}\,
\psi(2^{j_l} s_l-k_l)\,ds_l\\
&=& 2^{- j_l,H_l }\left\{\psi^{H_l}(2^{j_l}t_l-k_l)-\psi^{H_l}(-k_l)\right\}.
\end{eqnarray}
By inserting~(\ref{eq:KtWavDecomp}) into ~(\ref{Eq:Repstable}) for
every $l=1,\ldots, N$, we get that for any fixed  $t\in\R^N$, the
series (\ref{eq:wav-dec}) converges in probability to $X_0(t)$ and Proposition
\ref{thm:wav-dec} follows from Remark~\ref{rem:modif-lfss}.
\cqfd

From now on our goal will be to prove Proposition
\ref{prop:wav-holder}.  We need some preliminary results.

\noindent {\bf Proof of Proposition \ref{prop:wav-holder}.} \,
For the sake of simplicity we suppose that $N=2$. The
proof for the general case is similar. The space $C^\gamma(\K)$
is endowed with the norm
$$
\|f\|_\gamma=\sup_{x\in\K}|f(x)|+|f|_\gamma\quad\text{with}\quad|f|_\gamma
=\sup_{x\neq y\in\K}\frac{|f(x)-f(y)|}{\|x-y\|^\gamma} \; ,
$$
where $\|\cdot\|$ denotes the Euclidean norm in $\R^2$. For every $n\in\N$
we set $D_{n}^c=(\Z^2\times\Z^2)\setminus
D_{n,M}^2$.
Let us define $F_n (x,y)=F_n(x,y;\psi^{H_1};M,\phi,\delta,\beta,\eta)$
and $E(x,y)=E(x,y;\phi,\delta,\beta,\eta)$ by
$$
F_n (x,y)=A_n (x,y)+B_n(x,y) \; ,
$$
where $A_n(x,y)$ and $B_n(x,y)$ are defined in Lemma~\ref{lem:holder1}
in the Appendix, and
\begin{equation*}
E(x,y)=\sum_{(J,K)\in\Z^2}
  2^{-J\delta}\frac{|\phi(2^J x-K)-\phi(2^J
    y-K)|}{|x-y|^\beta}(3+|J|)^{1/\alpha+\eta}(3+|K|)^{1/\alpha+\eta}.
\end{equation*}
Using (\ref{eq:def-Un}), the triangle inequality and Lemma~\ref{lem:ub-eps},
one has for any $n,p\in\N$ and $s_1,s_2,t_1,t_2\in
[-M,M]$, denoting by $\prod$ the product over indices $l=1,2$,
\begin{equation*}
\begin{split}
&\frac{|U_{n+p,M}(s_1,s_2)-U_{n+p,M}(t_1,t_2)-U_{n,M}(s_1,s_2)
+U_{n,M}(t_1,t_2)|}{(|s_1-t_1|^2+|s_2-t_2|^2)^{\beta/2}}\\
&=\Big|
\sum_{(j,k) \in D_{n+p,M}^{2}\setminus D_{n,M}^{2}} 2^{-\l j,H \r}\, \epsilon_{j,k}
\frac{\prod \big(\psi^{H_l}(2^{j_l}s_l
-k_l)-\psi^{H_l}(-k_l)\big)-\prod\big(\psi^{H_l}(2^{j_l}t_l
-k_l)-\psi^{H_l}(-k_l)\big)}{(|s_1-t_1|^2+|s_2-t_2|^2)^{\beta/2}}
\Big| \\
&\le \sum_{(j,k) \in D_{n}^c} 2^{-\l j,H \r}\, |\epsilon_{j,k}|
\frac{\Big |\prod\big(\psi^{H_l}(2^{j_l}s_l-k_l)-\psi^{H_l}(-k_l)\big)-\prod\big(\psi^{H_l}(2^{j_l}t_l
-k_l)-\psi^{H_l}(-k_l)\big)\Big |}{(|s_1-t_1|^2+|s_2-t_2|^2)^{\beta/2}}\\
&\le C
\Bigg[\sum_{(j,k) \in D_{n}^c} 2^{-\l j,H \r}\,\prod (3+|j_l|)^{1/\a+\eta}
\prod (3+|k_l|)^{1/\a+\eta}\frac{\big|\psi^{H_1}(2^{j_1}s_1-k_1)-\psi^{H_1}(2^{j_1}t_1-k_1)\big|}{|s_1-t_1|^\beta} \\
&\hspace{5cm}\times\big|\psi^{H_2}(2^{j_2}s_2-k_2)-\psi^{H_2}(-k_2)\big| +\text{ similar term }\Bigg]\;,
\end{split}
\end{equation*}
where $C$ is a positive random variable. Observing that, for any non-negative array $(a_{j,k})_{(j,k)\in\Z^4}$,
$$
\sum_{(j,k) \in D_{n}^c}a_{j,k}\leq
\Big[\sum_{(j_1,k_1)\in\Z^2\setminus D_{n,M}}\quad \sum_{(j_2,k_2)\in\Z^2}
+\sum_{(j_1,k_1)\in\Z^2} \quad \sum_{(j_2,k_2)\in\Z^2 \setminus D_{n,M}}\Big]
a_{j,k}\;,
$$
we thus get
\begin{equation*}
\begin{split}
&\frac{|U_{n+p,M}(s_1,s_2)-U_{n+p,M}(t_1,t_2)-U_{n,M}(s_1,s_2)
+U_{n,M}(t_1,t_2)|}{(|s_1-t_1|^2+|s_2-t_2|^2)^{\beta/2}}\\
&\le C
\Big[F_n(s_1,t_1;\psi^{H_1};H_1,\beta,\eta)E(s_2,0;\psi^{H_2};H_2,0,\eta)
+E(s_1,t_1;\psi^{H_1};H_1,\beta,\eta)F_n(s_2,0;\psi^{H_2};H_2,0,\eta)\\
& \ \  +
F_n(s_2,t_2;\psi^{H_2};H_2,\beta,\eta)E(t_1,0;\psi^{H_1};H_1,0,\eta)
+E(s_2,t_2;\psi^{H_2};H_2,\beta,\eta)F_n(t_1,0;\psi^{H_1};H_1,0,\eta)\Big].
\end{split}
\end{equation*}
By Lemma~\ref{lem:holder1}, we have that $\sup_{x,y\in[-M,M]}F_n(x,y)\to 0$
as $n\to\infty$ and $\sup_{x,y\in[-M,M]}E(x,y)<\infty$; hence the last
display yields that $\sup_{p\geq0}|U_{n+p,M}-U_{n,M}|_\gamma\to0$ as
$n\to\infty$. Observing that $U_{n,M}$ vanishes on the axes, the same
result holds with $|\cdot|_\gamma$ replaced by $\|\cdot\|_\gamma$ and
Proposition~\ref{prop:wav-holder} is proved.
\cqfd

\begin{remark}
\label{rem:rem-holgauss} {\rm Proposition \ref{prop:wav-holder} is
much easier to prove  in the Gaussian case. Indeed, in this
case, using the fact that the $\epsilon_{j,k}$'s are {\em
independent} $\mathcal{N}(0,1)$ Gaussian random variables one can easily show
that the sequence $(U_{n,M})_{n\in\N}$ is weakly relatively compact in
the space $C(\K)$. We refer to the proof of Proposition 3 in
\cite{AX05} for more details.}
\end{remark}

From now on we will always identify the LFSS $X_0$ with its
random wavelet series representation (\ref{eq:wav-dec}).

\section{Uniform modulus of continuity and behavior as
$|t_\ell|\rightarrow 0$ or $\infty$}

The goal of this section is to prove Theorem \ref{thm:mod-cont}. An
immediate consequence of Proposition~\ref{prop:wav-holder} is that $X_0$
is locally $C^\gamma$ for any $\gamma\in(0,H_1-1/\a)$, almost surely.
Theorem~\ref{thm:mod-cont} completes this result by
providing a sharper estimate on the uniform modulus of continuity,
see~(\ref{eq:mod-cont}), and the behavior at infinity and around the
axes, see~(\ref{eq:ub-infty}). As in our note~\cite{ayache:roueff:xiao:2007a},
these results are obtained by using the
representation~(\ref{eq:wav-dec}). However, we improved the modulus
of continuity estimate by relying
on the independence present in the coefficients $\{\epsilon_{j,k},
\, (j,k) \in\Z^N \times\Z^N\}$, see
Lemma~\ref{lem:up-ps-indep}. If this independence is not taken into
account, an alternative result (i.e., Lemma~\ref{lem:ub-eps}) may be used,
resulting in a less precise estimate. The latter result
holds in a quite general framework since they can
be extended to a more general class of random wavelet series,
see Remark~\ref{rem:sub-infty} below.

\noindent {\bf Proof of Theorem \ref{thm:mod-cont}.}\, It follows from
(\ref{eq:wav-dec}), Corollary~\ref{cor:ub}
 and Lemma~\ref{lem:estim-ST} that for every
$\omega\in\Omega_{0}^*$ and every $s,t\in K$, the triangle
inequality implies
\begin{eqnarray}
\nonumber
\big|X_0(s,\omega)-X_0(t,\omega) \big| &\le& \sum_{n=1}^N
\big|X_0(t_1,\ldots, t_{n-1}, s_n,\ldots, s_N;\omega)-X_0(t_1,\ldots, t_n,s_{n+1},\ldots,
s_N;\omega) \big|\\
\nonumber
&\le& C_1(\omega) \sum_{n=1}^N\bigg(\prod_{\ell=1}^{n-1}
T_{H_\ell,1/\a,\eta}(t_\ell; \psi^{H_\ell})\bigg) \times
\bigg(\prod_{\ell=n+1}^N T_{H_\ell,1/\a,\eta}(s_\ell;\psi^{H_\ell})\bigg)\\
\label{eq:TSboundIncX}
&& \qquad \qquad  \times S_{H_n,1/\a,\eta}(t_n,s_n;\psi^{H_n})\\
\nonumber
&\le&  C_2(\omega) \sum_{n=1}^N |t_n-s_n|^{H_n-1/\alpha}\,
\big(1+\big|\log |t_n-s_n|\big| \big)^{2/\a+2\eta}.
\end{eqnarray}
This shows~(\ref{eq:mod-cont}).

Similarly, using~(\ref{eq:wav-dec}), Corollary~\ref{cor:ub}
and Lemma~\ref{lem:estim-ST}, we obtain, for every
$\omega\in\Omega_{0}^*$ and every $t\in \R$,
\begin{equation}
\label{eq:unifBoundX0}
\big|X_0(t,\omega)\big| \leq C_3(\omega) \prod_{\ell=1}^N
T_{H_\ell,1/\a,\eta}(t_\ell;\psi^{H_\ell}) \leq
C_4(\omega) \prod_{\ell=1}^N  \big(1+ \big| \log |t_\ell| \big|\big)^{1/\a+\eta}
|t_\ell|^{H_\ell}.
\end{equation}
The proof of Theorem \ref{thm:mod-cont} is finished.
\cqfd

\begin{remark}
\label{rem:sub-infty}
{\rm Clearly Proposition~\ref{prop:wav-holder} holds more generally for
any process $Y=\{Y(t),\, t\in\R^N\}$ having a wavelet series
representation of the form
$$
Y(t)=\sum_{(j,k)\in\Z^N\times\Z^N}c_{j,k}\lambda_{j,k}\prod_{l=1}^N
\left\{\phi_l(2^{j_l}t_l-k_l)-\phi_l(-k_l)\right\},
$$
where the $\phi_l$'s are well-localized functions,
$\{c_{j,k},\,j,k\in\Z^N\}$ is a sequence of complex-valued
coefficients satisfying $|c_{j,k}|\le c 2^{-\l j,H \r}$
($c>0$ being a constant) and  $\{\lambda_{j,k},\,j,k\in\Z^N\}$
is a sequence of complex-valued random
variables satisfying $\sup_{j,k}\E [|\lambda_{j,k}|^\nu ]<\infty$
for all $0<\nu<\alpha$.
We can also show that~(\ref{eq:ub-infty}) holds with probability 1 for
such a process $Y$. In contrast, for this more general class of process,
we cannot show~(\ref{eq:mod-cont}) but a less
precise estimate for the uniform modulus of continuity. Namely, as announced
in our note~\cite{ayache:roueff:xiao:2007a}, with
probability 1,
$$
\sup_{s,t\in {\mathcal K}}\frac{|X_0(s,\omega)-X_0(t,\omega)|}{\sum_{j=1}^N
  |s_j-t_j|^{H_j-1/\a-\eta}}<\infty \;
$$
for all compact sets ${\mathcal K} \subseteq \R^N$.}
\end{remark}

\section{Optimality of the modulus of continuity estimate}

The goal of this section is to prove Theorem \ref{thm:irreg}. For
every $n\in\{1,\ldots,N\}$ and $(j_n,k_n)\in\N\times\Z$, let
$G_{j_n,k_n}=\{G_{j_n,k_n} (\widehat{u}_n),\,\widehat{u}_n\in\R^{N-1}\}$
be the $\alpha$-stable field defined as the following wavelet transformation:
\begin{equation}
\label{eq:proj-X}
G_{j_n,k_n}(\widehat{u}_n)=2^{j_n(1+H_n)}\int_{\R}X_0(s_n,\widehat{u}_n)
\psi^{-H_n}(2^{j_n}s_n-k_n)\;ds_n,
\end{equation}
where the notation $(s_n,\widehat{u}_n)$ is introduced in Theorem~\ref{thm:irreg}.
By using (\ref{eq:ub-infty}) and the fact that the wavelet
$\psi^{-H_n}$ is well-localized, the process $\{G_{j_n,k_n}(u),\,u\in\R^{N-1}\}$ is
well-defined and its trajectories are continuous, almost surely. The proof of
Theorem \ref{thm:irreg} mainly relies on the following Lemmas
\ref{lem:contrad-holder} and \ref{lem:liminf-G}.
\begin{lemma}
\label{lem:contrad-holder} Let $\Omega_{0}^*$ be the event of
probability $1$ in Corollary~\ref{cor:ub} and let $n\in\{1,\ldots, N\}$.
Suppose that there exist $(u_{n},\widehat{u}_{n})\in\R^N$, $\rho>0$,
$\epsilon>0$ and $\omega\in\Omega_{0}^*$ such that
\begin{equation}
\label{eq:contrad-holder1}
\sup_{s_{n},t_{n}\in [u_{n}-\rho,u_{n}+\rho]}\frac{|X_0
(s_{n},\widehat{u}_{n},\omega)-X_0 (t_{n},\widehat{u}_{n},\omega)|}
{|s_{n}-t_{n}|^{H_{n}-1/\alpha} \big(1+\big|\log |s_{n}-t_{n}|\big| \big)^{-1/\a-\epsilon}}
<\infty.
\end{equation}
Then one has
\begin{equation}
\label{eq:contrad-holder2}
\limsup_{j_{n}\rightarrow\infty} \; (j_n 2^{-j_n})^{1/\alpha}
\max\big\{|G_{j_{n},k_{n}}(\widehat{u}_{n},\omega)|\,: \,k_{n} \in \Z,\,
|u_{n}-2^{-j_{n}}k_{n}|\le \rho/ 8\big\} =0.
\end{equation}
\end{lemma}

\begin{lemma}
\label{lem:liminf-G} Let $\Omega_{3}^*$ be the event of probability
$1$ defined as $\Omega_{3}^*=\Omega_{0}^*\cap \Omega_{2}^*$, where
$\Omega_{0}^*$ and $\Omega_{2}^*$ are respectively the events defined in
Corollary~\ref{cor:ub} and Lemma~\ref{lem:liminf-nu}. For all
$\omega\in\Omega_{3}^*$, $n\in\{1,\ldots,N\}$, all integers
$j_n\in\N$, real numbers $z_1< z_2$ and all  $0<\tau_1<\tau_2$,  one has
\begin{equation}
\label{eq:liminf-G} \liminf_{j_n\rightarrow\infty} \; (j_n2^{-j_n})^{1/\a}
\inf_{\widehat{u}_n\in [\tau_1,\tau_2]^{N-1}} \, \max\bigg\{
|G_{j_n,k_n}(\widehat{u}_n,\omega)|;\;k_n\in
[2^{j_n}z_1,2^{j_n}z_2]\cap\Z  \bigg\}>0.
\end{equation}
\end{lemma}

Before proving these lemmas, we show how they yield Theorem~\ref{thm:irreg}.

\noindent {\bf Proof of Theorem \ref{thm:irreg}.}\,
For the sake of simplicity we only consider the case where $\widehat{u}_n$
have positive and non-vanishing coordinates. The
general case is similar.
Suppose \textit{ad absurdum} that there exists
$\omega\in \Omega_{3}^*$ such that~(\ref{eq:irreg}) is not
satisfied. Then, for some $n\in\{1,\ldots ,N\}$, there exists
$\widehat{u}_{n}\in\R^{N-1}$ with positive and non vanishing
coordinates, some real number $u_n$, $\rho>0$ and $\epsilon>0$ arbitrary
small such that~(\ref{eq:contrad-holder1}) holds. By Lemma \ref{lem:contrad-holder},
this implies~(\ref{eq:contrad-holder2}). Then the conclusion of Lemma
\ref{lem:liminf-G} leads to a contradiction. This proves Theorem
\ref{thm:irreg}.
\cqfd

\noindent {\bf Proof of Lemma \ref{lem:contrad-holder}.}\,
Let $j_{n}\in\N$ and $k_{n}\in\Z$ be such that
\begin{equation}
\label{eq:cond1-jk}
|u_{n}-2^{-j_{n}}k_{n}|\le\rho/8.
\end{equation}
It follows from (\ref{eq:proj-X}) and (\ref{eq:moment-psiH}) that
$G_{j_{n},k_{n}}(\widehat{u}_{n},\omega)$ can be written as
\[
2^{j_{n}(1+H_{n})} \int_{\R} \Big ( X_0
(s_{n},\widehat{u}_{n},\omega)-X_0(2^{-j_{n}}k_{n},
\widehat{u}_{n},\omega)\Big ) \psi^{-H_{n}}
(2^{j_{n}}s_{n}-k_{n})\,ds_{n}.
\]
Hence, we have
\begin{eqnarray}
\nonumber
\big|G_{j_{n},k_{n}}(\widehat{u}_{n},\omega)\big| &\le&
2^{j_{n}(1+H_{n})}\int_{\R} \big| X_0(s_{n},
\widehat{u}_{n},
\omega)-X_0(2^{-j_{n}}k_{n},\widehat{u}_{n},\omega)
\big|\;\big|\psi^{-H_{n}}(2^{j_{n}} s_{n}- k_{n}) \big|\,ds_{n}\\
\label{eq:bound1-G}
&=&2^{j_{n}(1+H_{n})}\left\{ A_{j_{n},k_{n}}(\widehat{u}_{n},\omega)+
B_{j_{n},k_{n}}(\widehat{u}_{n},\omega)\right\}\;,
\end{eqnarray}
where
\begin{equation}
\label{eq:def-GA}
A_{j_{n},k_{n}}(\widehat{u}_{n},\omega)
=\int_{|s_{n}-u_{n}|\le \rho/2}
\big|X_0(s_{n},\widehat{u}_{n},\omega)-
X_0(2^{-j_{n}} k_{n},\widehat{u}_{n},\omega)\big|\; \big|\psi^{-H_{n}}(2^{j_{n}}
s_{n}-k_{n})\big|\,ds_{n}
\end{equation}
and
\begin{equation}
\label{eq:def-GB}
B_{j_{n},k_{n}}(\widehat{u}_{n},\omega)
=\int_{|s_{n}-u_{n}|> \rho/2}
\big|X_0(s_{n},\widehat{u}_{n},\omega)-X_0(2^{-j_{n}}k_{n},
\widehat{u}_{n},\omega)\big|\;\big|\psi^{-H_{n}}(2^{j_{n}}s_{n}-k_{n})\big|\,ds_{n}.
\end{equation}

Let us now give a suitable upper bound for $A_{j_{n},k_{n}}
(\widehat{u}_{n},\omega)$. It follows from (\ref{eq:def-GA})
and (\ref{eq:contrad-holder1}) that
$A_{j_{n},k_{n}}(\widehat{u}_{n},\omega)$ is at most
\begin{equation}
\label{eq:Abound1}
C_{5}(\omega)\int_{\R}
\big|s_{n}-2^{-j_{n}}k_{n}\big|^{H_{n}-1/\alpha}\,
\big(1+\big|\log |s_n-2^{-j_{n}}k_{n}| \big| \big)^{-1/\a-\epsilon}
\big|\psi^{-H_{n}}(2^{j_{n}}s_{n}-k_{n})\big|\,ds_{n}.
\end{equation}
We claim that
\begin{equation}
\label{eq:logPSI}
\sup_{j_n\geq1}
\int_{\R}\big|x\big|^{H_{n}-1/\alpha}\big( 1/j_n+\big| \log 2- (\log |x|)/j_n\big|
\big)^{-1/\a-\epsilon}\,\big|\psi^{-H_{n}}(x)\big|\,dx<\infty
\end{equation}
and differ its proof after we have shown (\ref{eq:contrad-holder2}).

By setting $x=2^{j_{n}}s_{n}-k_{n}$ in the integral in~(\ref{eq:Abound1})
and using~(\ref{eq:logPSI}),
one obtains, for all $j_n\geq1$ and $k_n\in\Z$ satisfying~(\ref{eq:cond1-jk}),
\begin{equation}
\label{eq:GB-bound1}
A_{j_{n},k_{n}}(\widehat{u}_{n},\omega) \le C_6(\omega)
2^{j_{n}(-1-H_n+1/\alpha)}j_n^{-1/\a-\epsilon} \; .
\end{equation}

In order to derive an upper bound for
$B_{j_{n},k_{n}}(\widehat{u}_{n},\omega)$, we use the fact
that $\psi^{-H_{n}}$ is a well-localized function and~(\ref{eq:cond1-jk})  to get
\begin{eqnarray*}
B_{j_{n},k_{n}}(\widehat{u}_{n},\omega) &\le& c\,
\int_{|s_{n}-u_{n}|>\rho/2}
\big|X_0(s_{n},\widehat{u}_{n},\omega)
-X_0 (2^{-j_{n}}k_{n},\widehat{u}_{n},\omega) \big|\;
\Big(1+|2^{j_{n}}s_{n}-k_{n}|\Big )^{-2}\,ds_{n}\\
& \le&  c \,\int_{|s_{n}-u_{n}|>\rho/2}
\big|X_0(s_{n},\widehat{u}_{n},\omega)
-X_0(2^{-j_{n}}k_{n},\widehat{u}_{n},\omega)\big|\\
&& \hspace{5cm} \times \left(1+2^{j_n}\Big (|s_{n}-u_{n}|
- |u_{n}-2^{-j_{n}}k_{n}|\Big)\right)^{-2}\,ds_{n}\\
&\le& c\, 2^{-2j_{n}}
\int_{|s_{n}-u_{n}|>\rho/2}
\big|X_0(s_{n},\widehat{u}_{n},\omega) - X_0 (2^{-j_{n}}
k_{n}, \widehat{u}_{n},\omega) \big|\,
\big|s_{n}-u_{n}\big|^{-2}\,ds_{n}.
\end{eqnarray*}
This last inequality, together with (\ref{eq:ub-infty}), implies
that, since $\omega\in\Omega_0^*$,
\[
B_{j_{n},k_{n}}(\widehat{u}_{n},\omega) \le C_7(\omega) \, 2^{-2j_{n}},
\]
where $C_7$ is a random variable that does not depend on the integers $j_n$ and $k_n$
satisfying~(\ref{eq:cond1-jk}).
Hence, putting together the last inequality, (\ref{eq:GB-bound1})
and (\ref{eq:bound1-G}) one obtains~(\ref{eq:contrad-holder2}).

Finally, to conclude the proof of the lemma, it remains to show~(\ref{eq:logPSI}).
We separate the integral in (\ref{eq:logPSI}) into two domains, $|x|>2^{j_n/2}$ and
$|x|\leq2^{j_n/2}$.
We bound $\big(1/j_n+\big| \log 2-(\log |x|)/j_n\big| \big)^{-1/\a-\epsilon}$ from above by
$j_n^{1/\a+\epsilon}$ on the first
domain, and by $\big((\log 2)/2 \big)^{-1/\a-\epsilon}$ on the second domain, yielding that the
integral in~(\ref{eq:logPSI}) is at most
$$
j_n^{1/\a+\epsilon}\int_{|x|>2^{j_n/2}}
\big|x\big|^{H_{n}-1/\alpha}\big|\psi^{-H_{n}}(x)\big|\,dx
+ \big((\log 2)/2 \big)^{-1/\a-\epsilon} \int_{\R}\big|x\big|^{H_{n}-1/\alpha}\,
\big|\psi^{-H_{n}}(x)\big|\,dx \; .
$$
Using that $H_{n}-1/\alpha\in(0,1)$ and that $\psi^{-H_{n}}$ is well localized,
we thus get~(\ref{eq:logPSI}).
\cqfd

In order to prove Lemma \ref{lem:liminf-G}, we first prove a weaker
result, namely the following lemma.

\begin{lemma}
\label{lem:liminf-nu}
There exists $\Omega_{2}^*$, an event of probability $1$, such
that for all $\omega\in \Omega_{2}^*$, $n \in \{1,\ldots,N\}$  and
real numbers $M>1$, $z_1<z_2$, $0<\tau_1<\tau_2$, one has
\begin{equation}
\label{eq:liminf-nu} \liminf_{j_n\rightarrow\infty}\,
(j_n2^{-j_n})^{1/\alpha} \nu (n,j_n;M;z_1,z_2;\tau_1,\tau_2;\omega)>0\;,
\end{equation}
where
\begin{equation}
\label{eq:def-nu}
\begin{split}
&\nu (n,j_n;M;z_1,z_2;\tau_1,\tau_2;\omega)\\
& =\min_{\widehat{k}_n\in
[M^{j_n}\tau_1,M^{j_n}\tau_2]^{N-1}\cap\Z^{N-1}}  \max\bigg\{
\big|G_{j_n,k_n }(M^{-j_n}\widehat{k}_n,\omega)\big|;\;k_n\in
\Big[2^{j_n}z_1 ,\, 2^{j_n}z_2 \Big] \cap \Z
\bigg\} \;.
\end{split}
\end{equation}

\end{lemma}

In order to prove Lemma \ref{lem:liminf-nu} we need to show that
the random variables $G_{j_n,k_n}(\widehat{u}_n)$ satisfy some nice
properties, namely the following proposition.

\begin{proposition}
\label{prop:proper-G}
Let $\widehat{u}_n\in \R^{N-1}$ be an arbitrary
fixed vector with non-vanishing coordinates. Then the following results hold:
\begin{enumerate}
\item[(a)] $\{G_{j_n,k_n}(\widehat{u}_n),\;(j_n,k_n)\in\N\times\Z\}$ is a
sequence of strictly $\a$-stable random
variables with identical scale parameters given by
\begin{equation}
\label{eq:sigma-hatUn}
\|G_{j_n,k_n}(\widehat{u}_n)\|_\alpha=\|\psi\|_{L^\alpha (\R)} \prod_{l\neq n}
\left\|(u_l-\cdot)_{+}^{H_l-1/\alpha}-(-\cdot)_{+}^{H_l-1/\alpha}\right\|_{L^\alpha (\R)}
\end{equation}
\item[(b)] Let $L>0$ be a constant such that the support of $\psi$ is included in $[-L,L]$.
Then for all integers $p>2L$ and $j_n\ge 0$, $\{G_{j_n,q_n
p}(\widehat{u}_n);\;q_n\in\Z\}$ is a sequence of independent random
variables.
\end{enumerate}
\end{proposition}

Proposition \ref{prop:proper-G} is in fact a straightforward
consequence of the following proposition and the fact that any two functions
$s_n \mapsto \psi (2^{j_n}s_n - q_np)$ with different values of $q_n$
have disjoint supports.

\begin{proposition}
\label{prop:rep-G}
For every vector $\widehat{u}_n$ with non-vanishing coordinates and for
every $(j_n,k_n)\in\N\times\Z$ one has almost surely
\begin{equation}
\label{eq:rep-G} G_{j_n,k_n}(\widehat{u}_n)=\int_{\R^N}\bigg
[2^{j_n/\alpha} \psi (2^{j_n}s_n-k_n) \prod_{l\neq n} \Big\{
(u_l-s_l)_{+}^{H_l-1/\alpha} - (-s_l)_{+}^{H_l-1/\alpha}
\Big\}\bigg]\; dZ_\alpha (s)\;.
\end{equation}
\end{proposition}

\noindent {\bf Proof of Proposition \ref{prop:rep-G}.}\, As
in~(\ref{eq:unifBoundX0}), we have
\begin{equation}
\label{eq:sub-infty}
\sup_{t\in\R^N}\frac{\sum_{(j,k)\in\Z^N\times\Z^N}2^{-\l j,H\r}
\,|\epsilon_{j,k}|\Big |\prod_{l=1}^N
\left\{\psi^{H_l}(2^{j_l}t_l-k_l)-\psi^{H_l}(-k_l)\right\}\Big
|}{\prod_{j=1}^N
|t_j|^{H_j} \big(1+\big|\log |t_j| \big| \big)^{1/\alpha+\eta}}<\infty.
\end{equation}
It follows
from Propositions \ref{thm:wav-dec} and \ref{prop:wav-holder}, (\ref{eq:sub-infty}),
the Dominated Convergence Theorem,~(\ref{eq:biorth-psiH}),~(\ref{eq:moment-psiH})
and~(\ref{eq:def-kappa}) that for any $\widehat{u}_n\in\R^{N-1}$
one has almost surely, for any increasing sequence $(D_m)_{m\in\N}$
of finite sets in $\Z\times\Z$ such that $\cup_mD_m=\Z\times\Z$,
\begin{equation}\label{eq:decomp1-G}
\begin{split}
G_{j_n,k_n}(\widehat{u}_n)&=\lim_{m\rightarrow\infty} \sum_{(j,k)\in D_m^{N}}
2^{j_n-\l\widehat{j}_n,\widehat{H}_n\r}\epsilon_{j,k}\\
& \qquad \qquad \qquad \qquad\times \int_{\R}\prod_{l=1}^{N}
\Big[\psi^{H_l}(2^{j_l}s_l-k_l)-\psi^{H_l}(-k_l)\Big]
\psi^{-H_n}(2^{j_n}s_n-k_n)\,ds_n\\
&= \lim_{m\rightarrow\infty} \sum_{(\widehat{j}_n, \widehat{k}_n)\in
D_m^{N-1}} \prod_{l\neq n} \kappa_{l,j_l,k_l} (u_l)\,
\epsilon_{(j_n, \widehat{j}_n); (k_n,\widehat{k}_n)},
\end{split}
\end{equation}
where $\kappa_{l,j_l,k_l} (u_l)$ is defined in (\ref{eq:def-kappa}).
On the other hand, it follows from (\ref{eq:KtWavDecomp}) that
\begin{equation}
\label{eq:decomp-kernelG}
\psi_{j_n,k_n}(s_n)\prod_{l\neq n}
\Big\{(u_l-s_l)_{+}^{H_l-1/\alpha}-(-s_l)_{+}^{H_l-1/\alpha}\Big\}
=\sum_{(\widehat{j}_n,\widehat{k}_n)\in \Z^{2(N-1)}}\prod_{l\neq n}
\kappa_{l,j_l,k_l}(u_l) \prod_{l=1}^N \psi_{j_l,k_l}(s_l),
\end{equation}
where for all fixed $\widehat{u}_n\in \R^{N-1}$ the convergence of
the series in the RHS (\ref{eq:decomp-kernelG}), as a function of
$s\in\R^N$, holds in $L^\alpha(\R^N)$. Next using
(\ref{eq:decomp-kernelG}) and (\ref{eq:def-epsilon}) one has, for
every fixed $\widehat{u}_n\in\R^{N-1}$,
\begin{multline}
\label{eq:decomp2-G}
\int_{\R^N}\bigg[\psi_{j_n,k_n}(s_n)\bigg (\prod_{l\neq
n}\Big\{(u_l-s_l)_{+}^{H_l-1/\alpha}-(-s_l)_{+}^{H_l-1/\alpha}\Big\}\bigg
)\bigg ]\;dZ_\alpha (s)\\
=\sum_{(\widehat{j}_n,\widehat{k}_n)\in\Z^{2(N-1)}}
\prod_{l\neq n}\kappa_{l,j_l,k_l}(u_l)\epsilon_{(j_n,\widehat{j}_n);(k_n,\widehat{k}_n)},
\end{multline}
where the convergence of the series holds in probability. Finally, putting together
(\ref{eq:decomp1-G}), (\ref{eq:decomp2-G}) and (\ref{eq:def-psijk}), one
obtains the proposition.
\cqfd

We are now in position to prove Lemma \ref{lem:liminf-nu}.

\noindent {\bf Proof of Lemma \ref{lem:liminf-nu}.}\, For any constants $M, c_1 >0$, $n \in
\{1, \ldots, N\}$, integer $j_n\ge 0$ and rational numbers
$r_1<r_2$, $0<\theta_1<\theta_2$ and $\zeta >0$, let $\Gamma (n,j_n)=\Gamma
(n,j_n;M, c_1;r_1,r_2;\theta_1,\theta_2;\zeta)$ be the event defined as
\begin{equation}
\label{eq:ev-gamma}
\Gamma (n,j_n;M,c_1;r_1,r_2;\theta_1,\theta_2)=\Big\{\omega\,:\,\nu
(n,j_n;M;r_1,r_2;\theta_1,\theta_2;\omega)\le (c_1\,j_n2^{-j_n})^{-1/\alpha} \Big \}.
\end{equation}
First we will show that, there exists $c_1$ large enough such that
\begin{equation}
\label{eq:boca1}
\sum_{j_n\in\N}\P\Big (\Gamma (n,j_n;M,c_1;r_1,r_2;\theta_1,\theta_2)\Big)<\infty.
\end{equation}
Using~(\ref{eq:tailprob}),~(\ref{eq:sigma-hatUn}) and that
$\big\|(u_l-\cdot)_{+}^{H_l-1/\alpha}-(-\cdot)_{+}^{H_l-1/\alpha}\big\|_{L^\alpha
(\R)}$ is increasing with $|u_l|$ and non-zero for $u_l\neq0$, we have
\begin{equation}
\label{eq:havytail}
c_2 := \min_{n=1,\dots,N}\;\;\inf_{t\geq1}\;\;\inf_{(j_n,k_n)\in\N\times\Z}\;\;
\inf_{\widehat{u}_n\in[\theta_1,\theta_2]^{N-1}}
t^\a P(|G_{j_n,k_n}(\widehat{u}_n)|>t) >0 \;.
\end{equation}
Observe finally that
\begin{multline*}
\nu(n,j_n;M;r_1,r_2;\theta_1,\theta_2;\omega)\\
 \geq\min_{\widehat{k}_n \in [M^{j_n}\theta_1,
M^{j_n}\theta_2]^{N-1}\cap\Z^{N-1}} \max\bigg\{ \big|G_{j_n,q_n
p}(M^{-j_n}\widehat{k}_n,\omega)\big|;\;q_n\in
\Big[\frac{2^{j_n}r_1}{p},\, \frac{2^{j_n}r_2}{p}\Big ]\cap\Z
\bigg\}.
\end{multline*}
It follows from Proposition \ref{prop:proper-G} and (\ref{eq:havytail})
and this inequality that
\begin{equation}
\label{eq:prgamma1} \begin{split}
\P\Big (\Gamma (n,j_n)\Big)
& \le \sum_{\widehat{k}_n} \prod_{q_n\in [
\frac{2^{j_n}r_{1}} p, \frac{2^{j_n}r_{2}} p]\cap\Z}\P\Big
(|G_{j_n,q_n
p}(M^{-j_n}\widehat{k}_n)|\le (c_1\,j_n2^{-j_n})^{-1/\alpha}\Big )\\
&\le c_3\, M^{(N-1)j_n}\Big (1-c_2\, j_n\,2^{-j_n}/c_1\Big)^{c_4 2^{j_n}},
\end{split}
\end{equation}
where the summation is taken over all $\widehat{k}_n \in [M^{j_n}\theta_1,
M^{j_n}\theta_2]^{N-1}\cap\Z^{N-1}$ and the constants $c_2$, $c_3$ and $c_4$
do not depend on $j_n$. Using the last inequality one can prove that
(\ref{eq:boca1}) holds by choosing $c_1>0$ large enough. Hence the
Borel-Cantelli Lemma implies that, for such a constant $c_1$,
$$
\P\bigg (\bigcup_{m\in\N}\bigcap_{j_n\ge m}
\Gamma^c (n,j_n;M,c_1;r_1,r_2;\theta_1,\theta_2) \bigg)=1,
$$
where $\Gamma^c(n,j_n;M,c_1,r_1,r_2;\theta_1,\theta_2)$ denotes the
complement event of $\Gamma(n,j_n;M,c_1;r_1,r_2;\theta_1,\theta_2)$.
But this implies that the event
$$
\left\{\omega:\liminf_{j_n\to\infty}
(j_n2^{-j_n})^{1/\alpha} \nu(n,j_n;M;r_1,r_2;\theta_1,\theta_2;\omega)>0\right\}
$$
has probability 1. Finally setting $\Omega_{2}^*$ as the intersection
of such sets over
$\Big\{(M;r_1,r_2;\theta_1,\theta_2)\in\Q^5\,:\,M>0,\,r_1<r_2 \mbox{ and }
0<\theta_1<\theta_2\Big\}$,
one obtains the lemma.
\cqfd

The following proposition will allow us to derive Lemma
\ref{lem:liminf-G} starting from Lemma \ref{lem:liminf-nu}. Roughly
speaking it means that the increments of the random field $\{G_{j_n,k_n}
(\widehat{u}_n),\widehat{u}_n \in [\tau_1,\tau_2]^{N-1}\}$ can be
bound uniformly in the indices $j_n$ and $k_n$.
\begin{proposition}
\label{prop:incre-G} Let $\Omega_0^*$ be the event of probability
$1$ that was introduced in Corollary~\ref{cor:ub}. Then for any
reals $z_1<z_2$, $0<\tau_1<\tau_2$ and $\eta >0$ arbitrarily small,
there exists an almost surely finite random variable $C_8>0$ such that
for every $n\in\{1,\ldots,N\}$, $j_n\in\N$, $k_n\in [2^{j_n}z_1,\,
2^{j_n}z_2]$, $\widehat{u}_n\in [\tau_1,\, \tau_2]^{N-1}$,
$\widehat{v}_n\in [\tau_1,\, \tau_2]^{N-1}$ and
$\omega\in\Omega_0^*$, one has
\begin{equation}
\label{eq:incre-G}
\big|G_{j_n,k_n}(\widehat{u}_n,\omega)-G_{j_n,k_n}(\widehat{v}_n,\omega)\big|
\le C_8(\omega)\, 2^{j_n H_n}\sum_{l\neq
n}|u_l-v_l|^{H_l-1/\alpha-\eta}.
\end{equation}
\end{proposition}

\noindent {\bf Proof.}\,  Lemma~\ref{lem:estim-ST} applied to~(\ref{eq:TSboundIncX})
shows that, for all $\omega\in\Omega_{0}^*$ and any $\eta>0$, there exists
$C(\omega)>0$ such that, for every
$n\in\{1,\ldots, N\}$, $s_n\in\R$, $\widehat{u}_n\in
[\tau_1,\tau_2]^{N-1}$ and $\widehat{v}_n\in [\tau_1,\tau_2]^{N-1}$,
\begin{equation}
\label{eq:modcont-infty}
\big|X_0(s_n,\widehat{u}_n,\omega)-
X_0(s_n,\widehat{v}_n,\omega)\big| \le C(\omega)\bigg (\sum_{l\neq
n}|u_l-v_l|^{H_l-1/\alpha-\eta}\bigg)\, |s_n|^{H_n}(1+|\log
(|s_n|)|)^{1/\alpha+\eta}.
\end{equation}
Let $\zeta>0$ be arbitrary small and consider the integral
$$
I(j_n,k_n)=2^{j_n}\int_{\R}(1+|s_n|)^{H_n +\zeta}\,
\big|\psi^{-H_n}(2^{j_n}s_n-k_n)\big|\,ds_n \;.
$$
By setting $x=2^{j_n}s_n-k_n$ we derive that
\begin{eqnarray}
\nonumber
\sup_{j_n\in\N}\max_{k_n\in [2^{j_n}z_1,\,2^{j_n}z_2]}I(j_n,k_n)
&=& \sup_{j_n\in\N}\max_{k_n\in [2^{j_n}z_1,\,2^{j_n}z_2]}
\int_{\R}\Big(1+2^{-j_n}|x+k_n|\Big)^{H_n+\zeta}\,
\big|\psi^{-H_n}(x)\big|\,dx\\
\label{eq:bound-intI}
&\le & \int_{\R}\Big(1+|x|+\max\{|z_1|,|z_2|\}\Big)^{H_n+\zeta}\,
\big|\psi^{-H_n}(x)\big|\,dx<\infty.
\end{eqnarray}
The inequality (\ref{eq:incre-G}) then follows from  (\ref{eq:proj-X}),
(\ref{eq:modcont-infty}) and (\ref{eq:bound-intI}).
\cqfd

We are now in position to prove Lemma \ref{lem:liminf-G}.

\noindent {\bf Proof of Lemma \ref{lem:liminf-G}.}\, We set
\begin{equation}
\label{eq:nu-tilde}
\widetilde{\nu}(n,j_n;z_1,z_2;\tau_1,\tau_2;\omega)
=\inf_{\widehat{u}_n\in [\tau_1,\tau_2]^{N-1}} \, \max\bigg\{
|G_{j_n,k_n}(\widehat{u}_n,\omega)|;\;k_n\in
[2^{j_n}z_1,2^{j_n}z_2]\cap\Z  \bigg\}.
\end{equation}
In view of Lemma \ref{lem:liminf-nu} it is sufficient to show that there
exists $\gamma>0$ small enough and $M>0$ such that, for all
$n\in\{1,\ldots,N\}$, $\omega\in\Omega_3$ and reals $z_1<z_2$,
$0<\tau_1<\tau_2$, one has
\begin{equation}
\label{eq:diff-nunutilde}
\lim_{j_n\rightarrow\infty}2^{-j_n(1/\alpha-\gamma)}
\Big|\widetilde{\nu}(n,j_n;M;z_1,z_2;\tau_1,\tau_2;\omega)
-\nu(n,j_n;z_1,z_2;\tau_1,\tau_2;\omega)\Big|=0.
\end{equation}
As the function $f_{j_n}(\cdot)= \max\Big\{
\big|G_{j_n,k_n}(\cdot,\omega) \big|; \;k_n\in [2^{j_n}z_1,\,
2^{j_n} z_2]\Big\}$ is continuous, there exists
$\widehat{u}_{n}^0(j_n)\in [\tau_1,\, \tau_2]^{N-1}$ such that
\begin{equation}
\label{eq:point-min} f_{j_n}(\widehat{u}_{n}^0(j_n))=\inf
\Big\{f_{j_n}(\widehat{u}_{n});\,\widehat{u}_n\in
[\tau_1,\tau_2]^{N-1}\Big\}.
\end{equation}
Moreover, when $j_n$ is big enough, one has for some $\widehat{k}_{n}^0(j_n)\in
[M^{j_n}\tau_1,\, M^{j_n}\tau_2]^{N-1}\cap\Z^{N-1}$,
\begin{equation}
\label{eq:approx-min}
\|M^{-j_n}\widehat{k}_{n}^0(j_n)-\widehat{u}_{n}^0(j_n)\|_\infty\le M^{-j_n} \;.
\end{equation}
Then it follows from
Proposition \ref{prop:incre-G} that there exists a constant $c_5>0$
(independent of $(j_n, k_n)$) such that the following inequality
holds
$$
\big|G_{j_n,k_n}(M^{-j_n}\widehat{k}_{n}^0(j_n),\omega)-G_{j_n,k_n}(\widehat{u}_{n}^0
(j_n),\omega) \big|\le c_5\, 2^{j_n H_n}M^{-j_n
(H_1-1/\alpha-\eta)}.
$$
The last inequality implies that
\begin{equation}
\label{eq:diff-max}
f_{j_n}(M^{-j_n}\widehat{k}_{n}^0(j_n)) \le
f_{j_n}(\widehat{u}_{n}^0(j_n))+c_5\, 2^{j_n H_n}M^{-j_n
(H_1-1/\alpha-\eta)}.
\end{equation}

By using (\ref{eq:point-min}) and (\ref{eq:diff-max}) one obtains
that
\begin{equation}\begin{split}
f_{j_n}(\widehat{u}_{n}^0(j_n))& \le \min
\Big\{f_{j_n}(M^{-j_n}\widehat{k}_{n});\,\widehat{k}_{n}\in
[M^{j_n}\tau_1,M^{j_n}\tau_2]^{n-1}\Big\}\\
&\le f_{j_n}(\widehat{u}_{n}^0(j_n))+c_5\, 2^{j_n H_n}M^{-j_n
(H_1-1/\alpha-\eta)}.
\end{split}
\end{equation}
Let us choose $M$ large enough so that
\begin{equation*}
\frac{H_N-1/\alpha}{H_1-1/\alpha}<\frac{\log M}{\log 2}.
\end{equation*}
and then, using (\ref{Eq:H}), we choose $\eta>0$ and $\gamma>0$
small enough so that
$2^{j_n H_n}M^{-j_n(H_1-1/\alpha-\eta)}=o(2^{-j_n(1/\a-\gamma)})$
as $j_n\to\infty$. Finally combining this with~(\ref{eq:diff-max}),
we obtain~(\ref{eq:diff-nunutilde}). This proves Lemma
\ref{lem:liminf-G}. \cqfd

\section{Proof of Theorem \ref{Th:dim}}

As usual, the proof of Theorem \ref{Th:dim} is divided into proving
the upper and lower bounds separately. The proofs of the lower
bounds rely on the standard capacity argument and the following
Lemma \ref{Lem:Compare}. However, the proofs of the upper bounds are
significantly different from that in \cite{AX05}, due to the fact
that both $\dim X\big([0, 1]^N\big)$ and $\dim {\rm Gr}X\big([0,
1]^N\big) $ are not determined by the exponent for the uniform
modulus of continuity of $X$. Our argument is based on the moment
method in \cite{LX94}. Combining this argument with the methods
in \cite{Xiao07}, we are able to determine the Hausdorff dimension
of the image $X(E)$ for all nonrandom Borel sets $E \subseteq (0, \infty)^N$.

We start by proving some results on the scale parameters of the increments
of real-valued LFSS $X_0$ between two points (i.e., $X_0(s) - X_0(t)$) and over
intervals; see Lemmas \ref{Lem:Compare} and \ref{Lem:IncRec} below.
Combining the latter with the maximal moment inequality due to Moricz
\cite{Moricz:1983}, we derive sharp upper bounds for the moments of the
supremum of $X_0$.

Lemma \ref{Lem:Compare} is an extension of Lemma 3.4 in \cite{AX05}
for fractional Brownian sheets. Since $d(s, t) :=\big\|X_0(s) - X_0(t)\big\|_\a $
can be used as a pseudometric for characterizing the regularity
properties of $X_0$ via metric entropy methods (cf. \cite[Chapter 12]{ST94}),
these results will be useful for studying other properties of LFSS $X$
as well.

\begin{lemma}\label{Lem:Compare}
For any constant $\ep> 0$, there exist positive and finite
constants $c_6$ and  $c_7$ such that for all $s, t
\in [\ep, 1]^N$,
\begin{equation}\label{Eq:IncVar}
c_6\, \sum_{\ell=1}^N \big|s_\ell - t_\ell\big|^{H_\ell} \le
\big\|X_0(s) - X_0(t)\big\|_\a \le c_7\, \sum_{\ell=1}^N
\big|s_\ell - t_\ell\big|^{H_\ell}.
\end{equation}
\end{lemma}
\noindent {\bf Proof.}\, To prove the upper bound in
(\ref{Eq:IncVar}), we use induction on $N$. When $N=1$, $X_0$ is an
$(H, \alpha)$-linear fractional stable motion and one can verify
directly that (\ref{Eq:IncVar}) holds as an equality. Suppose the
upper bound in (\ref{Eq:IncVar}) holds for any linear fractional
stable sheet with $n$ parameters. We now show that it holds for a
linear fractional stable sheet $X_0$ with $n+1$ parameters.

It follows from the representation (\ref{Eq:Repstable}) that, for
any $s, t \in [\ep, 1]^{n+1}$, $\big\|X_0(s) - X_0(t)\big\|_\a^\a $
is a constant multiple of the following integral:
\begin{equation*}
\begin{split}
&\int_{\R^{n+1}} \bigg| \prod_{\ell=1}^{n+1}
\Big\{(t_\ell-r_\ell)_+^{H_\ell- \frac 1 \a} -
(-r_\ell)_+^{H_\ell- \frac 1 \a} \Big\}-  \prod_{\ell=1}^{n+1}
\Big\{(s_\ell-r_\ell)_+^{H_\ell- \frac 1 \a} -
(-r_\ell)_+^{H_\ell- \frac 1 \a} \Big\} \bigg|^\a\, dr\\
& \le c\,  \int_{\R^{n}} \bigg| \prod_{\ell=1}^{n}
\Big\{(t_\ell-r_\ell)_+^{H_\ell- \frac 1 \a} -
(-r_\ell)_+^{H_\ell- \frac 1 \a} \Big\}-  \prod_{\ell=1}^{n}
\Big\{(s_\ell-r_\ell)_+^{H_\ell- \frac 1 \a} -
(-r_\ell)_+^{H_\ell- \frac 1 \a} \Big\} \bigg|^\a \; dr\\
&\qquad \qquad \times \int_{\R} \Big\{(t_{n+1}-r_{n+1})_+^{H_{n+1}- \frac 1
\a} -
(-r_{n+1})_+^{H_{n+1}- \frac 1 \a} \Big\}^\a \, dr_{n+1}\\
&+ c\,  \int_{\R^{n}} \bigg[ \prod_{\ell=1}^{n}
\Big\{(s_\ell-r_\ell)_+^{H_\ell- \frac 1 \a} -
(-r_\ell)_+^{H_\ell- \frac 1 \a} \Big\}\bigg]^\a  \; dr\\
&\qquad \qquad \times \int_\R \Big|(t_{n+1}-r_{n+1})_+^{H_{n+1}- \frac 1 \a} -
(s_{n+1}-r_{n+1})_+^{H_{n+1}- \frac 1 \a} \Big|^\a\, dr_{n+1}\\
&\leq c\,\left\{ \bigg(\sum_{\ell=1}^n |s_\ell -
t_\ell|^{H_\ell}\bigg)^\a + |t_{n+1} - s_{n+1}|^{H_{n+1}\a} \right\},
\end{split}
\end{equation*}
where, in deriving the last inequality, we have used the induction
hypothesis, the fact that the function $t \mapsto \int_\R
\{(t-r)_+^{H-1/\a}-(-r)_+^{H-1/\a}\}^\a dr$ is locally uniformly
bounded for any $H>1/\a$ and that, by a change of variable
$r_{n+1}=t_{n+1}+|t_{n+1}-s_{n+1}| u$, the last integral in the
previous display is less than $|t_{n+1} - s_{n+1}|^{\a H_{n+1}}$ up
to a multiplicative constant. Hence we have proved the upper bound
in~(\ref{Eq:IncVar}).

For proving the lower bound in (\ref{Eq:IncVar}), we define the
stable field $Y = \{Y(t), t \in \R_+^N\}$ by
\begin{equation}\label{Eq:Liouville}
Y(t) = \int_{[0,\, t]}h_{_H} (t, r)\, Z_\a(dr),
\end{equation}
where the function $h_{_H}(t, r)$ is defined in (\ref{Eq:functh}).
Then by  using (\ref{Eq:Repstable}) again we can write
\begin{equation}
\big\|X_0(s) - X_0(t)\big\|_\a \ge \big\|Y(t) - Y(s) \big\|_\a.
\end{equation}
To proceed, we use the same argument as in \cite[pp. 428--429]{AX05}
to decompose $Y$ as a sum of independent stable random fields. For
every $t \in [\ep, 1]^N$, we decompose the rectangle $[0, t]$ into
the following disjoint union:
\begin{equation}\label{Eq:Split0}
[0,\, t] = [0, \ep]^N \cup \bigcup_{j=1}^N R(t_j) \cup \Delta(\ep,
t),
\end{equation}
where $R(t_j) =\{ r \in [0,\, 1]^N: 0 \le r_i \le \ep \hbox{ if } i
\ne j,\, \ep < r_j \le t_j\}$ and $\Delta(\ep, t)$ can be written as
a union of $2^N - N -1$ sub-rectangles of $[0, t]$. It follows from
(\ref{Eq:Liouville}) and (\ref{Eq:Split0}) that for every $t \in
[\ep, 1]^N$,
\begin{eqnarray}\label{Eq:decomp}
Y(t) &=&  \int_{[0, \ep]^N} h_{_H} (t, r)\, Z_\a(dr) +
\sum_{j=1}^N \int_{R(t_j)}  h_{_H} (t, r)\, Z_\a(dr) +
\int_{\Delta(\ep,  t)}h_{_H} (t, r)\, Z_\a(dr) \nonumber\\
&:=& Y(\ep, t) + \sum_{j=1}^N Y_j(t) + Z(\ep, t).
\end{eqnarray}
Since the processes $\{Y(\ep, t),\,t\in\R^N\}$, $\{Y_j(t),\,
t\in\R^N\}$ ($1 \le j \le N$) and $\{Z(\ep, t),\,t\in\R^N\}$ are
defined by the stochastic integrals with respect to $Z_\a$ over
disjoint sets, they are independent. Only the $Y_j(t)$'s will be
useful for proving the lower bound in (\ref{Eq:IncVar}).

Now let $s, t \in [\ep,\, 1]^N$ and $ j \in \{1, \ldots, N\}$ be
fixed. Without loss of generality, we assume  $s_j \le t_j$. Then
\begin{equation}\label{Eq:Y3}
\big\|Y_j (t)- Y_j(s)\big\|_\a^\a = \int_{R(s_j)}\big(h_{_H}(t, r)
- h_{_H}(s, r)\big)^\a\, dr + \int_{R(s_j,  t_j)} h_{_H}^\a(t,
r)\, dr,
\end{equation}
where $R(s_j,  t_j) = \{ r \in [0, 1]^N: 0 \le r_i \le \ep \hbox{
if } i \ne j,\, s_j < r_j \le t_j\}$. By (\ref{Eq:Y3}) and some
elementary calculations we derive
\begin{equation} \label{Eq:Y4}
\begin{split}
\big\|Y_j (t)- Y_j(s)\big\|_\a^\a &\ge \int_{R(s_j, t_j)}
h_{_H}^\a(t, r) \, dr\\
&= \int_{[0, \ep]^{N-1}} \prod_{k\ne j} (t_k - r_k)^{\a H_k -1}
\int_{s_j}^{t_j} (t_j - r_j)^{\a H_j -1}\, dr\\ &\ge c\,
|t_j - s_j |^{\a H_j},
\end{split}
\end{equation}
where $c >0$ is a constant depending on $\ep$,
$\a$ and $H_k\, (1 \le k \le N)$ only.
The lower bound in (\ref{Eq:IncVar}) follows from (\ref{Eq:decomp}),
(\ref{Eq:Y3}) and (\ref{Eq:Y4}).
\cqfd

Our next lemma determines the scalar parameter of the increment of $X_0$
over any interval $[s, t] = \prod^N_{j=1}\, [s_j, t_j]$. Recall that
the increment of $X_0$ over $[s, t]$, denoted by $X_0([s,t])$, is defined as
\begin{equation}\label{Eq:X-inc}
X_0([s,t]) := \sum_{\delta \in \{0, 1\}^N} (-1)^{N - \sum_i \delta_i}
X_0(\l s_j + \delta_j (t_j - s_j)\r).
\end{equation}
This corresponds to the measure of the set $[s,t]$ by interpreting 
$X_0$ as a signed measure defined by
$X_0([0,t])=X_0(t)$ for all $t\in\R^N$ (convention: $[0, t_j] := [t_j, 0]$ 
if $t_j < 0$). 
It may be helpful to note that  for $N=2$, we have
\[
 X_0([s,t]) = X_0(t) - X_0((s_1, t_2)) - X_0((t_1, s_2)) + X_0(s).
\]
Similarly, we will denote the increment of the function $h_H(\cdot, r)$
over $[s, t]$ by $h_H([s, t], r)$.
\begin{lemma}\label{Lem:IncRec}
For any interval $[s, t] = \prod^N_{j=1}\, [s_j, t_j]$, we have
\begin{equation}\label{Eq:IncRec1}
\left\|X_0([s,t])\right\|_\a^\a = \prod_{j=1}^N (t_j - s_j)^{\a H_j}.
\end{equation}
\end{lemma}

\noindent{\bf Proof.}\, Since the kernel $h_H(t, s)$ in (\ref{Eq:functh}) is a tensor product,
it can be verified that
\begin{equation}\label{Eq:IncRec2}
\begin{split}
\left\|X_0([s,t])\right\|_\a^\a &= \int_{\R^N} \big| h_H([s, t], r)\big|^\a\, dr\\
&= \int_{\R^N}\kappa^\a \prod_{\ell=1}^N \Big|(t_\ell-r_\ell)_+^{H_\ell- \frac 1 \a}
- (s_\ell - r_\ell)_+^{H_\ell- \frac 1 \a} \Big|^\a\, dr\\
& =\prod_{j=1}^N (t_j - s_j)^{\a H_j}.
\end{split}
\end{equation}
This proves Lemma \ref{Lem:IncRec}. \cqfd

In order to estimate $\E\big[\sup_{t \in T} |X_0(t) - X_0(a)|\big]$
for all intervals $T = [a, b]\subseteq [\ep, 1]^N$, we will make use of
a general moment inequality of M\'oricz \cite{Moricz:1983} for the
maximum partial sums of multi-indexed random variables. This approach
has the advantage that it is applicable to non-stable random fields
as well. Another way for proving Lemma \ref{Lem:stablesuptail} below is to
establish sharp upper bounds for the tail probability
$\P\big[\sup_{t \in T} |X_0(t) - X_0(a)| > u\big]$ by modifying the
arguments in \cite{RosinskiSam93}.

First we adapt some notation from \cite{Moricz:1983} to our setting.
Let $\{\xi_k,\, k \in \N^N\}$ be a sequence of random variables. For
any $m \in \Z_+^N$ ($\Z_+$ is the set of nonnegative
integers) and $k \in \N^N$, let $R = R(m, k) = (m, m+k]\cap \Z_+^N$,
which will also be called a rectangle in $\Z_+^N$, and we denote
\begin{equation}\label{Eq:SR}
S(R) = S(m, k) = \sum_{p \in R}\xi_p \ \ \ \hbox{ and }\ \ \
M(R) = \max_{1 \le q \le k}\big|S(m, q)\big|.
\end{equation}
It can be verified that $M(R) \le \max_{Q \subseteq R} \big|S(Q)\big|
\le 2^N\, M(R)$, where the maximum is taken over all rectangles $Q
\subseteq R$. Let $f(R)$ be a nonnegative function of the rectangle $R$
with left-lower vertex in $\Z_+^N$. We call $f$ superadditive
if for every rectangle $R = R(m, k)$ the inequality
\begin{equation}\label{Eq:superadd}
f(R_{j 1}) + f(R_{j2}) \le f(R)
\end{equation}
holds for every $1 \le j \le N$ and $1 \le q_j < k_j$, where
\[
R_{j 1} = R\big((m_1, \ldots, m_N),\, (k_1, \ldots, k_{j-1}, q_j, k_{j+1},
\ldots, k_N)\big)
\]
and
\[
R_{j2} = R\big((m_1, \ldots, m_{j-1}, m_j+q_j, m_{j+1}, \ldots, m_N),\,
(k_1, \ldots, k_{j-1}, k_j-q_j, k_{j+1}, \ldots, k_N)\big).
\]
In other words, $R_{j1}\cup R_{j2}=R$ is a disjoint decomposition of $R$ by a
hyperplane which is perpendicular to the $j$th axis. Together with the
nonnegativity of $f$, (\ref{Eq:superadd}) implies that, for every fixed
$m \in \Z_+^N$, $f(R(m, k))$ is nondecreasing in each variable $k_j$
($1\le j \le N$).

The following moment inequality for the maximum $M(R)$ follows from
Corollary 1 in \cite{Moricz:1983}.

\begin{lemma}\label{Lem:Moricz}
Let $\beta > 1$ and $\gamma \ge 1$ be given constants. If there exists a
nonnegative and superadditive function $f(R)$ of the rectangle $R$ in $\Z_+^N$
such that $\E\big[|S(R)|^\gamma\big] \le f^{\beta}(R)$ for every $R$, then
\begin{equation}\label{Eq:Moricz1}
\E\big[M(R)^\gamma\big] \le \Big(\frac 5 2\Big)^N
\big(1 - 2^{(1- \beta)/\gamma}\big)^{N\gamma}\, f^\beta (R)
\end{equation}
for every rectangle $R$ in $\Z_+^N$.
\end{lemma}

It is useful to notice that the constant in \eqref{Eq:Moricz1} is independent of $R$.
Applying Lemma \ref{Lem:Moricz} to the linear fractional stable sheets, we obtain
\begin{lemma}\label{Lem:stablesuptail}
Let the assumption (\ref{Eq:H}) hold. Then there exists a positive
and finite constant $c_8$, independent of the skewness intensity $\beta(s)$, 
such that for all rectangles $T =
[a, b] \subseteq [\ep,1]^N$,
\begin{equation}\label{Eq:Moricz2}
\E\bigg(\sup_{t \in T} \big|X_0(t) - X_0(a)\big|\bigg) \le c_8\,
\sum_{j=1}^N (b_j - a_j)^{H_j}.
\end{equation}
\end{lemma}

\noindent {\bf Proof.}\, We prove this lemma by using induction on $N$.
In the case of $N=1$ it is well known (cf. \cite{LX94} or \cite{takashima:1989}) that
(\ref{Eq:Moricz2}) holds. Observe that the term $\delta=\l1\r$ in the sum appearing 
in~(\ref{Eq:X-inc}) is $X_0(t)$, and
since $\sum_{\delta \in \{0, 1\}^N} (-1)^{N - \sum_i \delta_i}=0$, we have
\begin{align}
\nonumber
X_0([s,t]) &= \sum_{\delta \in \{0, 1\}^N} (-1)^{N - \sum_i \delta_i}
\left(X_0(\l s_j + \delta_j (t_j - s_j)\r)-X_0(s)\right)\\
\label{Eq:X-inc2}
& = X_0(t)-X_0(s) + \sum_{\delta \in \{0, 1\}^N\setminus\l1\r} (-1)^{N - \sum_i \delta_i}
\left\{X_0(\l s_j + \delta_j (t_j - s_j)\r)-X_0(s)\right\} \; .
\end{align}
For every $\delta \in \{0, 1\}^N\setminus\l1\r$, there is some $n\in\{1,\dots,N\}$ 
such that $\delta_n=0$· Observing that, using the notation introduced in Theorem~\ref{thm:irreg}, 
$\widehat{u}_n \mapsto X_0(a_n,\widehat{u}_n)$ is an $(N-1,1)$-LFSS as defined by~(\ref{Eq:Repstable}) 
and~(\ref{Eq:functh}) but with $N$ replaced by $N-1$, $\kappa$ multiplied by a
constant only depending on $a_n$ and bounded independently of $a_n$ (since $a_n\in[\ep,1]$) 
and a modified skewness intensity $\widetilde{\beta}$.
Hence using the induction hypothesis, we have, for every $\delta \in \{0, 1\}^N\setminus\l1\r$,
\begin{equation}\label{Eq:allNinductionHyp}
\E\bigg(\sup_{t \in T} \big|X_0(\l a_j + \delta_j (t_j - a_j)\r)-X_0(a)\big|\bigg) \le c\,
\sum_{j=1}^{N} (b_j - a_j)^{H_j}\;.
\end{equation}
Applying~(\ref{Eq:X-inc2}) with $s=a$ and~(\ref{Eq:allNinductionHyp}), the 
bound~(\ref{Eq:Moricz2}) is implied by
\begin{equation}\label{Eq:N2c}
\E\Big(\sup_{t \in T} \big|X_0([a, t])\big|\Big) \le c\, \prod_{j=1}^N (b_j - a_j)^{H_j}\;,
\end{equation}
which we are now going to prove.
This is where Lemmas \ref{Lem:IncRec} and \ref{Lem:Moricz} will be applied.

For all $n\in\N$ we define a grid in $[a,b]$ with
mesh $2^{-n}$ by the collection of points
$$
\tau_{n}(p)= \l a_j + p_j(b_j-a_j)2^{-n}\r,\quad p\in R(0,\l 2^n\r)
=\{1,\dots,2^n\}^N\;.
$$
For each $p\in R(0,\l 2^n\r)$, we define the random variable $\xi_p$ to be
the increment of $X_0$ over the elementary rectangle with upper-right vertex 
$\tau_{n}(p)$, $[\tau_{n}(p-\l1\r),
\tau_{n}(p)]$. Interpreting $X_0$ as a signed measure, we get
that for any rectangle $R(m,k)\subseteq R(0,\l 2^n\r)$ for $m\neq 0$,
\begin{equation}\label{Eq:xisum0}
\sum_{p \in R(m,k)} \xi_{p} = X_0([\tau_{n}(m),\, \tau_{n}(m+k)]).
\end{equation}

We are now ready to prove \eqref{Eq:N2c}. By the continuity of the sample
function $X_0(t)$ and the monotone convergence theorem, since the set
$\cup_{n\geq1}\{\tau_n(p)\;:\;p\in R(0,\l2^n\r)\}$ is dense in
$[a,b]$, it is sufficient to show
that for all integers $n \ge 1$,
\begin{equation}\label{Eq:Moricz3}
\E\bigg(\max_{p\in R(0,\l2^n\r)} \big|X_0([a, \,\tau_{n}(p)])\big|\bigg)
\le c_9\,\prod_{j=1}^N (b_j - a_j)^{H_j},
\end{equation}
where $c_9 >0$ is a finite constant independent of $[a,b]\subseteq
[\varepsilon,1]^N$ and $n$.

It follows from Lemma \ref{lem:up-ps-indep} in the Appendix that for
any strictly $\a$-stable random variable $Z$ with scale parameter 1 and
every $0 <\gamma < \alpha$, we have $\E(|Z|^\gamma) \le c_{10}$, where
$c_{10}$ depends on  $\a$ and $\gamma$ only. This fact,~(\ref{Eq:xisum0})
and Lemma \ref{Lem:IncRec}
imply that for every $1 < \gamma < \a$ and every rectangle $R= R(m,k)
\subseteq R(0,\l 2^n\r)$,
\begin{equation}\label{Eq:Moricz4}
\begin{split}
\E \bigg(\bigg|\sum_{p \in R} \xi_p\bigg|^\gamma\bigg)
&= \E\bigg(\Big| X_0([\tau_{n}(m),\, \tau_{n}(m+k)])\Big|^\ga \bigg)\\
&\le c_{10}\, \bigg[\prod_{j=1}^N \bigg(\frac{k_j (b_j - a_j)}
{2^n}\bigg)^{H_j}\bigg]^\ga\\
&\le \bigg[c_{11}\, \prod_{j=1}^N \bigg(\frac{k_j (b_j - a_j)}
{2^n}\bigg)^{H_j/H_1}\bigg]^{H_1\ga},
\end{split}
\end{equation}
where $c_{11} = c_{10}^{1/(H_1\ga)}$.
For every rectangle $R = R(m,k)$ included in $R(0,\l2^n\r)$,
let
\[
f(R) = c_{11}\, \prod_{j=1}^N \bigg(\frac{k_j (b_j - a_j)} {2^n}\bigg)^{H_j/H_1}.
\]
Note that, under assumption (\ref{Eq:H}), we have $\a \in (1,
2)$, $H_1 \a > 1$ and $H_j \ge H_1$ for $j=1,\dots, N$. Hence
the inequality $x^{H_j/H_1} + y^{H_j/H_1} \le (x+y)^{H_j/H_1}$ for
all $x, y > 0$ implies that $f$ is superadditive.

We take $\gamma \in (1, \a)$ such that $\beta =\ga H_1 >1$ and apply
Lemma \ref{Lem:Moricz} to derive
\begin{equation}\label{Eq:Moricz5}
\begin{split}
\E \bigg(\sup_{k \in R(0, \l 2^n \r) } \bigg|\sum_{p \in R(0,k)} \xi_p
\bigg|^\gamma\bigg)
&\le c_{12}\, \bigg[\prod_{j=1}^N (b_j - a_j)^{H_j/H_1}\bigg]^{H_1\ga}\\
&= c_{12}\, \bigg[\prod_{j=1}^N (b_j - a_j)^{H_j}\bigg]^{\ga},
\end{split}
\end{equation}
where $c_{12} >0$ is a finite constant independent of $[a, b]$ and $n$.
This proves \eqref{Eq:N2c} and thus Lemma \ref{Lem:stablesuptail}. \cqfd

We now proceed to prove Theorem \ref{Th:dim}.

\noindent{{\bf Proof of Theorem \ref{Th:dim}}.}\,  We only
prove~(\ref{Eq:range}), which is done by modifying the proof of
Theorem 4 in~\cite{AX05} and by making use of Lemmas
\ref{Lem:Compare} and \ref{Lem:stablesuptail}. The formula
(\ref{Eq:graph}) can be proven using similar arguments and we leave
it to the interested reader.

First we prove the lower bound in (\ref{Eq:range}). Let $\ep \in
(0,\, 1)$ be given and let $I = [\ep,\, 1]^N$. We will prove that
for every $0 < \gamma < \min \{ d, \ \sum_{\ell=1}^N \frac1 {H_\ell}
\}$, $\dim X(I) \ge \gamma$ almost surely. By Frostman's theorem 
(see \textit{e.g.}~\cite{falconer90} pages 64 and 65),
it
is sufficient to show that we have
\begin{equation}\label{Eq:Frostman}
\E\int_{I} \int_I \frac 1 {\|X(s) - X(t) \|^{\gamma}}\, ds dt <
\infty,
\end{equation}
where $\|\cdot\|$ denotes the Euclidean norm in $\R^d$.

It is known that for any $d$-dimensional distribution function $F$
in $\R^d$ with characteristic function $\varphi$ and any $\gamma
>0 $, we have
\begin{equation}\label{Eq:gammafunct}
2^{\gamma/2 -1}\Gamma\Big(\frac{\gamma} 2\Big) \int_{\R^d}
\|x\|^{-\gamma}\, F(dx) = (2\pi)^{-d/2} \int_0^{+\infty}
u^{\gamma-1} du \int_{\R^d} \exp\bigg(- \frac{\|x\|^2} 2\bigg)
\varphi(ux)\, dx.
\end{equation}
This equality can be verified by replacing $\varphi$ in the right
side of (\ref{Eq:gammafunct}) by its expression as a Fourier
integral and then performing a routine calculation. Applying
(\ref{Eq:gammafunct}) to the distribution of the stable random
variable $ \xi = \big(X(s) - X(t)\big)/\|X(s) - X(t) \|_\a$ and
using Fubini's theorem, we obtain
\begin{equation}\label{Eq:deno}
\begin{split}
\E\big(\|\xi\|^{-\gamma}\big) &\le c_{14}\, \int_{\R^d} \exp
\bigg(- \frac{\|x\|^2} 2\bigg)\, dx \int_0^\infty u^{\gamma -
1}\, \exp \Big(- c_{15} |u|^\a \|x\|^\a \Big)\, du \\
&= c_{16}\, \int_{\R^d} \exp \bigg(- \frac{\|x\|^2} 2
\bigg)\, \|x\|^{-\gamma} dx < \infty,
\end{split}
\end{equation}
where the last integral is convergent because $\gamma < d$.
Combining (\ref{Eq:deno}) with Lemma \ref{Lem:Compare} yields
\begin{equation}\label{Eq:Frostman2}
\E\int_{I} \int_I \frac 1 {\|X(s) - X(t) \|^{\gamma}}\, ds dt \le
\int_I \int_I \frac 1 {\big(\sum_{\ell=1}^N |s_\ell -
t_\ell|^{H_\ell} \big)^\gamma}\, ds dt < \infty,
\end{equation}
where the finiteness of the last integral is proved in \cite[p.
432]{AX05}. This proves (\ref{Eq:Frostman}) and hence the lower
bound in (\ref{Eq:range}).

To prove the upper bound in (\ref{Eq:range}), we use the covering
argument in \cite{LX94} and \cite{AX05}.
Since clearly $\dim X\big([0, 1]^N\big) \le d$ a.s. and Hausdorff dimension
is $\sigma$-stable \cite{falconer90}, it is sufficient to show that
for every $\ep \in (0, 1)$,
\begin{equation}\label{Eq:rangeupper}
\dim X\big([\ep, 1]^N\big) \le \sum_{j=1}^N \frac1 {H_j} \qquad {\rm
a.s.}
\end{equation}
This will be done by using a covering argument.

For any integer $n \ge 2$, we divide $[\ep, 1]^N$ into $m_n$
sub-rectangles $\{R_{n, i}\}$ with sides parallel to the axes and
side-lengths $n^{-1/H_j}$ ($j =1, \ldots, N$), respectively. Then
\begin{equation}\label{Eq:mn}
 m_n \le c\, n^{\sum_{j=1}^N \frac 1{H_j}}
\end{equation}
and $X\big([\ep, 1]^N\big)$ can be covered by $X(R_{n, i})$ ($1 \le
i \le m_n$). Denote the lower-left vertex of $R_{n, i}$ by $a_{n,
i}$.  Note that the image $X(R_{n, i})$ is contained in a rectangle
in $\R^d$ with sides parallel to the axes and side lengths at most
$2 \sup_{s \in R_{n, i}} \big|X_k(s) - X_k(a_{n, i})\big|$ ($ k =
1,\ldots, d$), respectively. Hence each $X(R_{n, i})$ can be
covered by at most
$$  \prod_{k=1}^d \bigg[\frac{2 \sup_{s \in R_{n, i}} \big|X_k(s) - X_k(a_{n,
i})\big|} {n^{-1}} + 1 \bigg]
$$
cubes of side-lengths $n^{-1}$. In this way, we have obtained a
$(\sqrt{d}\,n^{-1})$-covering for $X\big([\ep, 1]^N\big)$.

By Lemma \ref{Lem:stablesuptail}, we derive that for every $1 \le
i \le m_n$ and $1 \le k \le d$,
\begin{equation}\label{Eq:diam2}
\E\bigg( \sup_{s \in R_{n, i}} \big|X_k(s) - X_k(a_{n, i})\big|
\Big) \le c\, n^{-1}.
\end{equation}

It follows from (\ref{Eq:mn}), (\ref{Eq:diam2}) and the
independence of $X_1, \ldots, X_d$ that for any $\ga >
\sum_{j=1}^N \frac 1 {H_j}$, we have
\begin{equation}\label{Eq:diam3}
\begin{split}
\E \bigg\{&\sum_{i=1}^{m_n}  \prod_{k=1}^d \bigg[\frac{2\, \sup_{s
\in R_{n, i}} \big|X_k(s) - X_k(a_{n, i})\big|} {n^{-1}} + 1
\bigg]\,
\Big(\sqrt{d}\, n^{-1}\Big)^{\ga}\bigg\}\\
& \le c\, n^{\sum_{j=1}^N \frac 1 {H_j}} \, n^{-\ga} \to
0 \qquad \hbox{ as }\ n \to \infty.
\end{split}
\end{equation}
This and Fatou's lemma imply that $\dim X\big([\ep, 1]^N\big) \le \ga$
almost surely. By letting $\ga \downarrow \sum_{j=1}^N \frac 1
{H_j}$ along rational numbers, we derive (\ref{Eq:rangeupper}). This
completes the proof of Theorem \ref{Th:dim}. \cqfd

The above method can be extended to determine the Hausdorff dimension
of the image $X(E)$ for every nonrandom Borel set $E \subseteq (0, \infty)^N$,
thus extending the results in Wu and Xiao
~\cite{wu:xiao:2007} and Xiao \cite{Xiao07} for anisotropic Gaussian random
fields to $(N, d)$-LFSS.

For this purpose, let us first recall from \cite{Xiao07} the definition
of a Hausdorff-type dimension which is more convenient to capture the
anisotropic nature of $X$.

For a fixed $(H_1, \ldots, H_N) \in (0, 1)^N$, let $\rho$ be the metric
on $\R^N$ defined by
\begin{equation}\label{Def:metric}
\rho(s, t) = \sum_{j=1}^N |s_j - t_j|^{H_j},\qquad \forall\, s, t
\in \R^N.
\end{equation}
For any $\beta > 0$ and $E \subseteq \R^N$,
define the $\beta$-dimensional Hausdorff measure [in the metric $\rho$]
of $E$ by
\begin{equation}\label{Eq:rhomeas}
{\mathcal H}^\beta_\rho(E) = \lim_{\de \to 0}
\inf\bigg\{\sum_{n=1}^\infty (2 r_n)^\beta:\, E \subseteq
\bigcup_{n=1}^\infty B_\rho(r_n),\, r_n \le \de\bigg\},
\end{equation}
where $B_\rho(r)$ denotes a closed (or open) ball of radius $r$
in the metric space $(\R^N, \rho)$. Then ${\mathcal H}^\beta_\rho$
is a metric outer measure and all Borel sets are
${\mathcal H}^\beta_\rho$-measurable. The
corresponding Hausdorff dimension of $E$ is defined by
\begin{equation}\label{Eq:rhodim}
\dimr E = \inf\big\{\beta> 0:\, {\mathcal H}^\beta_\rho(E) =
0\big\}.
\end{equation}
We refer to \cite{Xiao07} for more information on the history and basic
properties of ${\mathcal H}^\beta_\rho$ and $\dimr$.

\begin{theorem} \label{Th:dimXE}
Let the assumption (\ref{Eq:H}) hold. Then, for every nonrandom Borel
set $E \subseteq (0, \infty)^N$,
\begin{equation}\label{Eq:rangeXE}
\dim X(E) = \min\big\{ d; \ \dimr E \big\}\qquad \hbox{ a.s.}
\end{equation}
\end{theorem}

\noindent {\bf Proof.}\, The proof is a modification of those of Theorem
\ref{Th:dim} above and Theorem 6.11 in \cite{Xiao07}. For any $\gamma >
\dimr E$, there is a covering $\{B_\rho(r_n), n\ge 1\}$ of $E$
such that $\sum_{n=1}^\infty (2r_n)^\gamma \le 1$. Note that $X(E)
\subseteq \cup_{n=1}^\infty X \big(B_\rho(r_n)\big)$ and we can cover each
$X \big(B_\rho(r_n)\big)$ as in the proof of Theorem
\ref{Th:dim}. The same argument shows that $\dim X(E) \le
\gamma$ almost surely, which yields the desired upper bound for $\dim X(E)$.

By using the Frostman lemma for  ${\mathcal H}^\beta_\rho$ (Lemma 6.10 in
\cite{Xiao07}) and the capacity argument in the proof of Theorem
\ref{Th:dim}, one can show  $\dim X(E) \ge \min\big\{ d; \ \dimr E \big\}$
almost surely. We omit the details.
\cqfd

\appendix
\renewcommand{\theequation}{\thesection.\arabic{equation}}

\section{Technical lemmas}
The following lemma allows to
control the growth of an arbitrary sequence of strictly $\a$-stable
random variables having the same scale parameter.
\begin{lemma}
\label{lem:ub-eps} Let $\{\epsilon_\lambda,\,\lambda\in\Z^d\}$ be an
arbitrary sequence  of strictly $\a$-stable random variables having
the same scale parameter. Then, there exists
 an event $\Omega_1^*$ of probability $1$, such that for any $\eta>0$ and any
$\omega\in\Omega_1^*$,
\begin{equation}
\label{eq:ub-eps}
|\epsilon_\lambda (\omega)|\le C(\omega)\prod_{l=1}^d
(3+|\lambda_l|)^{1/\alpha+\eta},
\end{equation}
where $C>0$ is an almost surely finite random variable, only depending on
$\eta$.
\end{lemma}
\noindent {\bf Proof}\, This lemma simply follows from the  fact
that for any $\nu\in ((1/\alpha+\eta)^{-1},\alpha)$ one has
$$
\E\left (\sup_{\lambda\in\Z^d}\frac{|\epsilon_\lambda |^\nu}{\prod_{j=1}^d
    (3+|\lambda_j|)^{\nu (1/\alpha +\eta)}}\right )\le
c \sum_{\lambda\in\Z^d}\prod_{j=1}^d (3+|\lambda_j|)^{-\nu (1/\alpha
    +\eta)}<\infty.
$$
\cqfd

\begin{lemma}
\label{lem:up-ps-indep}
Let $\a\in(0,2)$. There exists a constant $c_{17}$ depending only on $\a$
such that for any strictly $\a$-stable random variable $Z$ with scale
parameter $\|Z\|_\a>0$ and skewness parameter
$\beta\in[-1,1]$ and all $t\geq\|Z\|_\a$,
\begin{equation}
  \label{eq:tailprob}
  c_{17}^{-1} \,\|Z\|_\a^\a\,t^{-\a} \leq \P(|Z|>t)\leq c_{17} \,\|Z\|_\a^\a\,t^{-\a}.
\end{equation}
Let $N\geq1$. Suppose now that $\{Z_{j,k},\,j_\ell\geq1,\;k_\ell\geq2\;\;\text{for}\;\;\ell=1,\dots,N\}$
is a sequence of strictly
$\a$-stable random variables such that
\begin{enumerate}[(i)]
\item For all $j\in(\N\setminus\{0\})^N$, $\{Z_{j,k},\,k_\ell\geq2\;\;\text{for}\;\;\ell=1,\dots,N\}$ are
    independent;
\item For all $j\in(\N\setminus\{0\})^N$ and $k\in(\N\setminus\{0,1\})^N$,
$\|Z_{j,k}\|_\a\leq1$.
\end{enumerate}
Then, with probability 1, one has, for any $\gamma>0$,
\begin{equation}
  \label{eq:supprobalphasemiindep}
\sup\left\{|Z_{j,k}|\prod_{\ell=1}^Nj_\ell^{-1/\a-\gamma}
k_\ell^{-1/\a}\log^{-1/\a-\gamma} k_\ell\;:\;j_\ell\geq1,
\;k_\ell\geq2\;\;\text{for}\;\;\ell=1,\dots,N\right\}<\infty \; .
\end{equation}
\end{lemma}
\noindent {\bf Proof}\,
Relation~(\ref{eq:tailprob}) follows from Property~1.2.15 in \cite{ST94}.
Let us now show~(\ref{eq:supprobalphasemiindep}) for $N=1$, the proof for
$N>1$ is similar. By using~(\ref{eq:tailprob}), we
obtain, for all $j\geq1$ and $n\geq1$,
$$
\P\big(\max\{|Z_{j,2}|,\dots,|Z_{j,n}|\}> u_{j,n}\big)\leq 1-
(1- c_{17} \,u_{j,n}^{-\a})^n,
$$
where $u_{j,n}=j^{1/\a+\gamma}n^{1/\a}\log^{1/\a+\gamma}n$.
Defining $n_m=[\exp(m)]$, we obtain
$$
\E\left[\sum_{j\geq1}\sum_{m\geq1}\1_{\max\{|Z_{j,2}|,\dots,|Z_{j,n_m}|\}
> u_{j,n_m}}\right]=
\sum_{j\geq1}\sum_{m\geq1}\P\Big(\max\{|Z_{j,2}|,\dots,|Z_{j,n_m}|\}
> u_{j,n_m}\Big)<\infty \;.
$$
Thus the random variable $\sum_{j\geq1}\sum_{m\geq1}\1_{\max\{|Z_{j,2}|,
\dots,|Z_{j,n_m}|\}> u_{j,n_m}}$ is a.s. finite. As a consequence there
exists an a.s. finite positive random variable $C$ such that
$$
\max\{|Z_{j,2}|,\dots,|Z_{j,n_m}|\}\leq C\,u_{j,n_m} \quad\text{for all}
\quad j\geq1,\; m\geq 1 \; .
$$
Let $m(k)$ be the unique integer satisfying $n_{m(k)}\leq k < n_{m(k)+1}$.
Thus for all $j\geq1,k\geq2$, we have
$$
|Z_{j,k}| \leq C\,u_{j,n_{m(k)+1}} = C\,j^{1/\a+\gamma}\;n_{m(k)+1}^{1/\a}\;
\log^{1/\a+\gamma}(n_{m(k)+1}) \;,
$$
Observe now that we have,  for all $k\geq2$,
$$
n_{m(k)+1}\leq \exp(m(k)+1) \leq \rme\,(n_{m(k)}+1) \leq \rme\; (k+1)\;.
$$
Relation~(\ref{eq:supprobalphasemiindep}) follows from the last two displays.
\cqfd

\begin{lemma}
\label{lem:shiftBis}
For any $\gamma\in[0,1)$ and $\eta\geq0$, there exists a constant $c>0$ such that,
for all $u\in\R$,
$$
\sum_{k\in\Z}(2+|u-k|)^{-2}(1+|k|)^\gamma \log^{\eta}(2+|k|)\leq c\,(1+|u|)^{\gamma}
\, \log^{\eta}(2+|u|) \; .
$$
\end{lemma}
\noindent {\bf Proof}\,
Put $k'=[u]-k$, where $[u]$ is the integer part of $u$. Hence
\begin{align*}
\sum_{k\in\Z}\frac{(1+|k|)^{\gamma}\log^{\eta}(2+|k|)}{(2+|u-k|)^{2}}
&=\sum_{k'\in\Z}(2+|u-[u]+k'|)^{-2}(1+|[u]-k'|)^{\gamma}\log^{\eta}(2+|[u]-k'|)\\
&\leq\sum_{k'\in\Z}(1+|k'|)^{-2} (2+|u|+|k'|)^{\gamma}\log^{\eta}(2+|u|+|k'|)\;.
\end{align*}
The result then follows by observing that
$(2+|u|+|k'|)^{\gamma}
\leq (1+|u|)^{\gamma}(2+|k'|)^{\gamma}$, $\log^{\eta}(2+|u|+|k'|)\leq c\,
\log^{\eta}(2+|u|)\log^{\eta}(2+|k'|)\}$ and $\gamma-2<-1$.
\cqfd

\begin{lemma}
\label{lem:geolog}
Let $\theta\neq0$ and $\gamma\in\R$. Set
$c:=\sum_{n\geq0}2^{-|\theta|n}(1+n)^{|\gamma|} <\infty$.
Then for any $n_0< n_1$ in $\{0,\pm1,\pm2,\dots,\pm\infty\}$,
\begin{equation}
\label{eq:geolog}
\sum_{n=n_0}^{n_1} 2^{n\theta}(1+|n|)^\gamma \le c\,
\begin{cases}
2^{n_0 \theta}(1+|n_0|)^\gamma & \text{if $\theta<0$}\\
2^{n_1 \theta}(1+|n_1|)^\gamma & \text{if $\theta>0$}.
\end{cases}
\end{equation}
\end{lemma}

\noindent {\bf Proof.}\,Take \textit{e.g.} $\theta<0$ and write
$$
\sum_{n=n_0}^{n_1} 2^{n\theta}(1+|n|)^\gamma \leq 2^{n_0 \theta}(1+|n_0|)^\gamma
\sum_{m\geq0}2^{m\theta}\left(\frac{1+|m+n_0|}{1+|n_0|}\right)^\gamma \;.
$$
Now observe that
$$
\frac1{1+|m|}\leq\frac{1+|m+n_0|}{1+|n_0|} \leq 1+|m|
$$
so that $\sup_{n_0}\sum_{m\geq0}2^{m\theta}\left(\frac{1+|m+n_0|}{1+|n_0|}\right)^\gamma
<\infty$ for any $\gamma \in \R$.\cqfd

\begin{lemma}
\label{lem:holder1}
For any $M>0$, $\eta >0$ small enough, $\delta\in (1/\alpha+\eta,\,1)$,
$\beta\in [0,\, \delta-1/\alpha-\eta)$, any well-localized function $\phi$
and $x,y\in\R$, let $A_n (x,y):= A_n
(x,y;M,\phi,\delta ,\beta,\eta)$ be the quantity defined as
\begin{eqnarray}
\label{eq:def-A}
A_n (x,y)=\sum_{|J|\le n}\sum_{|K|>
  M2^{n+1}}2^{-J\delta}\frac{|\phi(2^J x-K)-\phi(2^J y-K)|}
  {|x-y|^\beta}(3+|J|)^{1/\alpha+\eta}(3+|K|)^{1/\alpha+\eta}
\end{eqnarray}
and let $B_n (x,y) := B_n (x,y;\phi;\delta ,\beta,\eta)$ be the quantity
defined as
\begin{eqnarray}
\label{eq:def-B}
B_n (x,y)=\sum_{|J|\ge n+1}\sum_{K\in\Z}
  2^{-J\delta}\frac{|\phi(2^J x-K)-\phi(2^J y-K)|}{|x-y|^\beta}
  (3+|J|)^{1/\alpha+\eta}(3+|K|)^{1/\alpha+\eta},
\end{eqnarray}
with the convention that $A_n (x,x)=B_n (x,x)=0$ for any $x\in\R$.
These quantities converge to $0$, uniformly in $x,y\in [-M,M]$,
as $n$ goes to infinity.
\end{lemma}

\noindent {\bf Proof.} Let $x,y\in[-M,M]$ and $J_0\ge -\log_2(2M)$
be the unique integer such that
\begin{equation}
\label{eq:def-holder-J0}
2^{-J_0-1}<|x-y|\le 2^{-J_0}.
\end{equation}
Let us first prove that $A_n (x,y)$ converges to $0$, uniformly in
$x,y$ as $n$ goes to infinity. From now on we suppose that
$J$ is an arbitrary integer satisfying $|J|\le n$. We need to derive
suitable  upper bounds for the quantity
\begin{equation}
\label{eq:def-AJ}
A_{n}^{(J)}(x,y)=\sum_{|K|>M2^{n+1}}\frac{|\phi (2^J x-K)-\phi (2^J y-K)|}
{|x-y|^\beta}
(3+|K|)^{1/\alpha +\eta}.
\end{equation}
For this purpose, we consider two cases $J\le J_0$ and  $J\ge J_0+1$
separately. First we suppose that
\begin{equation}
\label{eq:JlessJ0}
J\le J_0.
\end{equation}
Using the Mean Value Theorem, (\ref{eq:wl-frac}), (\ref{eq:def-holder-J0}) and
(\ref{eq:JlessJ0}) one obtains that
\begin{equation*}
\begin{split}
|\phi(2^J x-K)-\phi(2^J y-K)|&\le c\, 2^{J}|x-y|\sup_{u\in I}(3+|u|)^{-2}\\
& \le c\, 2^J |x-y|(2+|2^J x-K|)^{-2},
\end{split}
\end{equation*}
where $I$ denotes the compact interval  with end-points $2^J x-K$ and $2^J
y-K$, whose length is at most 1 by~(\ref{eq:def-holder-J0}) and~(\ref{eq:JlessJ0}).
Next the last inequality and~(\ref{eq:def-AJ}) entail that
\begin{equation}
\label{eq:AJbound1}
A_{n}^{(J)}(x,y)\le c\, 2^J |x-y|^{1-\beta}\sum_{|K|>
  M2^{n+1}}\frac{(3+|K|)^{1/\alpha+\eta}}{(2+|2^J x-K|)^2}.
\end{equation}
On the other hand, using that $|x|\le M$ and $|J|\le n$, for all
$|K|>M2^{n+1}$, one gets
\begin{equation}
\label{eq:AJbound2}
\frac{(3+|K|)^{1/\alpha+\eta}}{(2+|2^J x-K|)^2}\le
\frac{(3+|K|)^{1/\alpha+\eta}}{(2+|K|-M2^n)^2}\le
c\, \big(1+|K| \big)^{-(2-1/\alpha-\eta)}.
\end{equation}
Putting together (\ref{eq:AJbound1}), (\ref{eq:AJbound2})
and~(\ref{eq:def-holder-J0}), one obtains that
\begin{equation}
\label{eq:AJbound3}
A_{n}^{(J)}(x,y)\le c\, 2^{J_0(\beta-1)}
2^{J-n(1-1/\alpha-\eta)}.
\end{equation}
Let us now study the second  case where
\begin{equation}
\label{eq:JmoreJ0}
J\ge J_0+1.
\end{equation}
It follows from (\ref{eq:def-holder-J0}), (\ref{eq:JmoreJ0}) and (\ref{eq:def-AJ})
that
\begin{equation}
\label{eq:AJbound4}
A_{n}^{(J)}(x,y)\le 2^{J\beta}\sum_{|K|>M2^{n+1}}\Big\{|\phi(2^J
x-K)|+|\phi(2^J x-K)|\Big\}(3+|K|)^{1/\alpha+\eta}.
\end{equation}
On the other hand, using (\ref{eq:wl-frac}) and the fact that $|J|\le n$ one
has, for any real $u\in [-M, M]$ and any $K\in\Z$ satisfying
$|K|>M2^{n+1}$,
\begin{equation}
\label{eq:AJbound5}
\big|\phi(2^J u-K)\big|\le c\, (3+|2^J u-K|)^{-2}\le c\, (3+|K|-M2^n)^{-2}\le
c_{18}\, (3+|K|)^{-2}.
\end{equation}
Combining (\ref{eq:AJbound4}) with (\ref{eq:AJbound5}) one gets that
\begin{equation}
\label{eq:AJbound6}
A_{n}^{(J)}(x,y) \le c_{19}\, 2^{J\beta-n(1-1/\alpha-\eta)}\;.
\end{equation}
It follows from (\ref{eq:def-A}),~(\ref{eq:def-AJ}), (\ref{eq:AJbound3})
and (\ref{eq:AJbound6}) that
\begin{equation*}
\begin{split}
A_n(x,y)&  \le c\, 2^{-n(1-1/\alpha-\eta)}
\left[2^{J_0(\beta-1)}\sum_{J=-\infty}^{J_0}2^{J(1-\delta)}(3+|J|)^{1/\alpha+\eta}
+\sum_{J=J_0+1}^\infty 2^{J(\beta-\delta)}(3+|J|)^{1/\alpha+\eta}\right]\\
& \leq c\,2^{-n(1-1/\alpha-\eta)}2^{J_0(\beta-\delta)}(3+|J_0|)^{1/\alpha+\eta}\\
&\leq c_{20}\, 2^{-n(1-1/\alpha-\eta)}\,,
\end{split}
\end{equation*}
where we used Lemma~\ref{lem:geolog} to bound the series and then the fact
$2^{-J_0}\leq 2M$ (see~(\ref{eq:def-holder-J0})).
Since $c_{20}$ does not depend on $(x,y)$, the last inequality proves that $A_n (x,y)$
converges to $0$, uniformly in $x,y\in[-M,M]$ as $n$ goes to infinity.

Let us now prove that $B_n(x,y)$ converges
to $0$, uniformly in $x,y$ as $n$ goes to infinity. In all the sequel
$J$ denotes an arbitrary integer satisfying $|J|\ge n+1$. First, we derive a suitable
upper bound for the quantity
\begin{equation}
\label{eq:def-BJ}
B_{n}^{(J)}(x,y)=\sum_{K\in\Z}\frac{|\phi (2^J x-K)-\phi (2^J y-K)|}{|x-y|^\beta}
(3+|K|)^{1/\alpha +\eta}.
\end{equation}
As above, we distinguish two cases: $J\le J_0$ and $J\ge J_0+1$. First we
suppose that (\ref{eq:JlessJ0}) is verified. As in~(\ref{eq:AJbound1}), we have
$B_{n}^{(J)}(x,y)\le c\, 2^J
|x-y|^{1-\beta}\sum_{K\in\Z}(3+|K|)^{1/\alpha+\eta}(2+|2^J x-K|)^{-2}$.
Next, using~(\ref{eq:def-holder-J0}) and Lemma \ref{lem:shiftBis} and the fact
that $|x|\le M$, one obtains that
\begin{equation}
\label{eq:BJbound1}
B_{n}^{(J)}(x,y)\le c\, 2^{J+J_0(\beta-1)} (1+2^J)^{1/\alpha+\eta}.
\end{equation}
Now let us suppose that (\ref{eq:JmoreJ0}) is verified. By using this
relation, (\ref{eq:def-holder-J0}), the triangle inequality,
(\ref{eq:wl-frac}), Lemma \ref{lem:shiftBis} and the fact that $x,y\in [-M,M]$,
one gets
\begin{eqnarray}
\label{eq:BJbound2}
\nonumber
B_{n}^{(J)}(x,y)&\le & 2^{J\beta}\sum_{K\in\Z}\Big\{|\phi(2^J x-K)|+|\phi(2^J
y-K)|\Big\}(3+|K|)^{1/\alpha+\eta}\\\nonumber
&\le & c\, 2^{J\beta}\sum_{K\in\Z}\Big\{(3+|2^J x-K|)^{-2}+(3+|2^J
y-K|)^{-2}\Big\}(3+|K|)^{1/\alpha+\eta}\\\nonumber
&\le & c\, 2^{J\beta}\Big\{(1+2^J |x|)^{1/\alpha+\eta}
+(1+2^J |y|)^{1/\alpha+\eta}\Big\}\\
&\le & c\, 2^{J(\beta+1/\alpha+\eta)}.
\end{eqnarray}

Since $2^{-J_0}\le M$, for all $n\geq\log_2(2M)$, we have $-n\leq J_0$, and thus, by~(\ref{eq:BJbound1}),
\begin{eqnarray}
\nonumber
\sum_{J\leq-n}2^{-J\delta}(3+|J|)^{1/\a+\eta}B_{n}^{(J)}(x,y)
& \leq&  c\, 2^{J_0(\beta-1)}\,\sum_{J\leq-n}2^{J(1-\delta)}(3+|J|)^{1/\a+\eta}\\
\label{eq:Bbound1}
&\leq &c\,  2^{n(\delta-1)}(1+n)^{1/\a+\eta}\;,
\end{eqnarray}
where we used Lemma~\ref{lem:geolog} and  $2^{-J_0}\le M$. Applying Lemma~ \ref{lem:geolog}
with~(\ref{eq:BJbound1}) and~(\ref{eq:BJbound2}) yields
\begin{equation}
  \label{eq:Bbound2}
  \sum_{J\in\Z}2^{-J\delta}(3+|J|)^{1/\a+\eta}B_{n}^{(J)}(x,y)
\le
c\, 2^{J_0(\beta+1/a+\eta-\delta)}(3+|J_0|)^{1/\a+\eta} \; ,
\end{equation}
and for any $n\geq J_0$,
\begin{equation}
  \label{eq:Bbound3}
\sum_{J\geq n}2^{-J\delta}(3+|J|)^{1/\a+\eta}B_{n}^{(J)}(x,y)
\le
c\, 2^{n(\beta+1/a+\eta-\delta)}(3+n)^{1/\a+\eta} \; .
\end{equation}
Since $\beta+1/a+\eta-\delta<0$, the function
$t\mapsto2^{t(\beta+1/a+\eta-\delta)}(3+t)^{1/\a+\eta}$ is decreasing for $t$
large enough, and hence for $n$ large enough, either $n\geq J_0$ and we may
apply~(\ref{eq:Bbound3}), or $n\leq J_0$ and we may apply~(\ref{eq:Bbound2})
whose right-hand side is smaller than the right-hand side of~(\ref{eq:Bbound3}).
Hence~(\ref{eq:Bbound3}) holds for all $n$ large
enough independently of $J_0$. This, with~(\ref{eq:Bbound1}), shows that $B_n(x,y)$
converges uniformly in $x,y$, as $n$ goes to infinity.
\cqfd

\begin{lemma}
\label{lem:estim-ST} Let $\phi$ be a well-localized  function i.e. a
function satisfying the condition~(\ref{eq:wl-frac}). For any
$\delta\in (0,1)$, $\gamma\in(0,\delta)$ and $\eta\geq0$, define
\begin{equation}
\label{eq:def-S} S_{\delta,\gamma,\eta}(x,y;\phi)=
\sum_{(J,K)\in\Z^2}2^{-J\delta}|\phi (2^J x-K)-\phi
(2^Jy-K)|(3+|J|)^{\gamma+\eta}(3+|K|)^\gamma\log^{\gamma+\eta}(2+|K|)
\end{equation}
and
\begin{equation}
\label{eq:def-T} T_{\delta,\gamma,\eta}(x;\phi)=\sum_{(J,K)\in\Z^2}
2^{-J\delta}\, \big|\phi (2^J x-K)-\phi(-K)\big|
(3+|J|)^{\gamma+\eta}(3+|K|)^\gamma\log^{\gamma+\eta}(2+|K|).
\end{equation}
Then, there exists a constant $c>0$, only depending on $\delta$,
$\gamma$ and $\phi$, such that the inequalities
\begin{multline}
\label{eq:estim-ST1} S_{\delta,\gamma,\eta}(x,y;\phi)\le
c\,|y-x|^{\delta-\gamma}\big[|y-x|^{\gamma}+|x|^\gamma+|y|^\gamma\big] \\
\times \big(1+\big|\log |y-x| \big| \big)^{2\gamma+2\eta} \,\big\{\log^{\gamma+\eta}(2+|x|)+
\log^{\gamma+\eta}(2+|y|)\big\}
\end{multline}
and
\begin{equation}
\label{eq:estim-ST2}
T_{\delta,\gamma,\eta}(x;\phi)\le c \big (1+\big|\log |x| \big|
\big)^{\gamma+\eta}\, |x|^\delta
\end{equation}
hold for all $x,y\in\R$ (with the convention that $0^a\times\log^b 0=0$
for all $a,b>0$).
\end{lemma}
\noindent {\bf Proof.}\, We only
prove~(\ref{eq:estim-ST1}), the proof of~(\ref{eq:estim-ST2}) is similar.
By~(\ref{eq:wl-frac}), there is a constant $c>0$ such that, for all
$J,K\in\Z$ and $x,y\in\R$,
\begin{equation}
\label{eq:bound-phi1}
\big|\phi (2^J x-K)-\phi (2^Jy-K)\big| \leq c\,  \left\{(2+|2^Jx-K|)^{-2}
+ (2+|2^Jy-K|)^{-2}\right\}.
\end{equation}
The quantity $|\phi (2^J x-K)-\phi (2^Jy-K)|$ can be bounded more sharply
when the condition $2^J |x-y|\le 1$ holds, namely by using (\ref{eq:wl-frac})
and the Mean Value Theorem one obtains that
\begin{eqnarray}
\label{eq:bound-phi2}
\nonumber
|\phi (2^J x-K)-\phi (2^Jy-K)| &\leq & c\, 2^{J}|x-y|\sup_{u\in I}(3+|2^J u-K|)^{-2}\\
&\leq & c \, 2^J|x-y| (2+|2^Jx-K|)^{-2}  \; ,
\end{eqnarray}
where $I$ denotes the compact interval whose end-points are $x$ and
$y$. From now on we will assume that $x\ne y$ (Relation~(\ref{eq:estim-ST1})
is trivial otherwise) and let $J_0\in\Z$  be the unique integer satisfying
\begin{equation}
\label{eq:encad-xy}
1/2 <2^{J_0}|y-x|\leq 1.
\end{equation}
The inequalities (\ref{eq:bound-phi1}) and (\ref{eq:bound-phi2}) entail that
\begin{equation}
\label{eq:bound-S}
S_{\delta,\gamma,\eta}(x,y;\phi)\le c\, \big (A_{J_0}|x-y|+B_{J_0}\big ),
\end{equation}
where
\begin{equation*}
A_{J_0}=\sum_{J\le J_0}\sum_{K\in\Z}
2^{J(1-\delta)}(2+|2^J x-K|)^{-2}(3+|J|)^{\gamma+\eta}(3+|K|)^\gamma
\log^{\gamma+\eta}(2+|K|)
\end{equation*}
and
\begin{equation*}
B_{J_0}=\sum_{J>J_0}\sum_{K\in\Z}
2^{-J\delta}\Big\{(2+|2^J x-K|)^{-2} +(2+|2^J y-K|)^{-2}\Big
\}(3+|J|)^{\gamma+\eta}(3+|K|)^\gamma\log^{\gamma+\eta}(2+|K|).
\end{equation*}
Lemma~\ref{lem:shiftBis} and Lemma~\ref{lem:geolog} yield
$$
A_{J_0}\leq c\, 2^{J_0(1-\delta)}\,(1+|x|^\gamma2^{J_0\gamma})\,(1+|J_0|)^{2\gamma+2\eta}\,
\log^{\gamma+\eta}(2+|x|)
$$
and, since $\gamma-\delta<0$,
$$
B_{J_0} \leq c\,
2^{-J_0\delta}\,(1+(|x|^\gamma+|y|^\gamma)2^{J_0\gamma})\,(1+|J_0|)^{2\gamma+2\eta}\,
\{\log^{\gamma+\eta}(2+|x|)+\log^{\gamma+\eta}(2+|y|)\}
\; .
$$
Inserting these two bounds into (\ref{eq:bound-S}) and using~(\ref{eq:encad-xy}),
we get~(\ref{eq:estim-ST1}) and the proof
is finished.
\cqfd

\bigskip

\bigskip

\bibliographystyle{plainnat}
\bibliography{lfss}

\end{document}